\documentclass[12pt]{article}
\usepackage{a4,verbatim,amsmath,amssymb}
\usepackage[isolatin]{inputenc}
\usepackage{epsfig}
\input math1.def
\newtheorem{theor}{Theorem}[section]

\newtheorem{lemma}[theor]{Lemma}

\newtheorem{prop}[theor]{Proposition}

\numberwithin{equation}{section}
\def\proof{\goodbreak\noindent{\sc Proof. }\nobreak}

\def\endproof{\par\nobreak\hbox to \hsize{\hfil\vrule width 5pt height
5pt}\goodbreak\vskip 3pt}

\begin{document}
\title{A Decomposition for Hardy Martingales. Part  II.}
     
\author{Paul F. X. M\"uller\thanks{Supported
by the Austrian Science foundation (FWF) Pr.Nr. FWFP23987000. 
Participant of the NSF sponsored  {\it Workshop in
Analysis and Probability} at Texas A $\&$ M University 2012.}  }
\date{September $15$,  2012}
\maketitle
\begin{abstract} 
This paper continues \cite{pfxm12} on Davis and Garsia  
Inequalities (DGI).
We prove DGI for dyadic 
perturbations of  Hardy martingales,
and apply them  to estimate the    $L^1  $
distance of  
a dyadic martingale on $\bT^\bN $ to the class of 
Hardy martingales. We revisit Bourgain's embedding of $L^1$ 
into the quotient space  $L^1/ H^1_0 . $  The Appendix
reviews well known estimates 
 on cosine-martingales complementary to DGI \cite{b1}.

\end{abstract}
\tableofcontents
\section{Introduction}\label{intro}
The introduction is divided into three separate parts. 
We first list preliminary material, standard notations and tools 
employed throughout this paper. Then we survey the Davis and Garsia
inequalities for Hardy martingales and their dyadic perturbations.  
We discuss the essential steps of the proof, and point out the role 
DGI are playing in Bourgain's embedding of $L^1 $ into $L^1 / H^1_0 . $ 
\subsubsection*{Preliminaries, Notation, Conventions}
Let  $\bT=  \{ e^{i\theta} : \theta \in
[0, 2\pi [ \} $ be  the torus 
equipped with the  normalized angular measure.
 Let $\bT^\bN = \{(x_i)_{i=1}^\infty  \}$ its  countable product
 equipped with its product Haar measure $\bP .$ We let 
$\bE $ denote expectation  with respect to   $\bP .$ 

\paragraph{Martingales on  $\bT^\bN$. }
Denote by $\cF_n $
the  sigma-algebra on $\bT^\bN$ generated by the 
cylinder sets 
$ \{(A_1, \dots, A_n , \bT^\bN )\},$
where $A_i,\, i \le n $ are measurable subsets 
of $\bT . $ 
Thus $(\bT^\bN, (\cF_k) , \bP) $ becomes a filtered probability space.
We let $\bE_{n}$ denote the conditional expectation with
respect to the $\s-$algebra  $\cF_n .$
Let  $F = (F_k) $ be  an $L^1(\bT^\bN)-$bounded  martingale.
Conditioned on  $\cF_{k-1}$  the martingale difference $\Delta F_k =
F_k- F_{k-1}$  defines an element in  $L_0^1(\bT) ,$ the Lebesgue space of integrable,
 functions with vanishing mean.
\paragraph{Dyadic martingales. }
The dyadic sigma-algebra on $ \bT^\bN $ is  defined with 
 Rademacher functions. For $ x = ( x _k ) \in  \bT^\bN $
define $\cos_k ( x ) =  \Re x_k  $ and

$$ 
 \s_ k ( x ) = {\rm sign} (\cos_k ( x )).
$$
We let  $\cD$ be   sigma- algebra 
 generated by 
$ \{\s_k , k \in \bN \} . $
Let   $D = (D_k) $ be  an integrable $ (\cF_k)$ martingale in
$\bT^\bN .$
If each  $ D_k $ is  measurable with respect to  $\cD ,$
then we  say that  $D = (D_k) $ is a dyadic martingale on $\bT^\bN .$
Define the non negative kernel  
$$B (w,z) = \prod_{k = 1 }^\infty  ( 1 +  \s_k( w) \s_k ( z )) ,
\quad  w  , z  \in \bT^\bN . $$
Conditional expectation $ \bE _\cD $ is the integral 
operator with
kernel $B ,$ 
$$   \bE _\cD  G ( w ) =    \bE_z  ( B (w,z)  G(z)) , \quad \quad   G \in L^1 ( \bT^\bN ).$$ 

\paragraph{Hardy martingales. }
Let  $H^1_0 (\bT ) \sb  L^1_0 (\bT ) $   consist of those integrable
functions for which the harmonic extension to the
unit disk is analytic.   And put $H^\infty_0 (\bT ) = H^1_0 (\bT ) \cap
L^\infty  (\bT ) . $ See  \cite{jg81}.

An  $L^1(\bT^\bN ) $ bounded  $(\cF_n)$ martingale $F = (F_k) $
is called a Hardy martingale if conditioned on $\cF_{n-1}$ 
the martingale difference  
$$ \Delta F _n  = F_n - F_{n-1}  $$
defines an element in $H^1_0 (\bT ).  $  See \cite{g1}, \cite{gar2}.
\paragraph{Sine and Cosine martingales.}
A cosine martingale    $U = ( U_k ) $  on  $\bT^\bN $ 
is defined by the relation 
\begin{equation}\label{14aug12hh}
 \Delta U_k (  x , \overline{y} ) 
=  \Delta U_k (  x, y ) , \quad\quad x \in \bT^{k-1} , \, \quad y \in \bT . 
\end{equation} 
A sine-martingale    $V = ( V_k ) $  on  $\bT^\bN $ 
is defined by 
\begin{equation}\label{14aug12h}
 \Delta V_k (  x , \overline{y} ) 
= - \Delta V_k (  x, y ) , \quad\quad x \in \bT^{k-1} , \, \quad y \in \bT . 
\end{equation}

\paragraph{Classical martingale spaces.}
Let  
$ ( \O , ( \cF_n) , \bP ) $  be a  filtered probability spaces. 
Corresponding to the fixed filtration we define 
the well known  spaces of martingales  $ H^1 ,$ $ \cP $ and  $ \cA $
by specifying their norms. 
Let   $G = (G_k) $ be an integrable $ ( \cF_n) $ martingale.
Define the  previsible norm $\cP ,$ 
\begin{equation}\label{15augmt2}
\| G\|_\cP  =\bE ( \sum_{k = 1 }^n \bE_{k-1} |\Delta G_k |^2 )^{1/2},
\end{equation}
 the $H^1 $ norm, and  the absolutely summing norm $ \cA , $ 
\begin{equation}\label{15augmt3}
\| G\|_{H^1} =
\bE ( \sum_{k = 1 }^n  |\Delta G_k |^2)^{1/2}
\quad {\rm and} \quad 
\| G\|_{\cA} =
\bE ( \sum_{k = 1 }^n  |\Delta G_k |).
\end{equation}
We refer to $$( \sum_{k = 1 }^n \bE_{k-1} |\Delta G_k |^2 )^{1/2}$$
as the conditional square function of $G.$  
See \cite{sia}  for the  classical
inclusions,
$$ \cA \sbe H^1 \sbe  L^1, 
\quad {\rm and}  \quad \cP \sbe  H^1 . $$  
The last inclusion is the content of  
 the  Burkholder-Gundy inequality \cite{sia},  
\begin{equation}\label{150812y}
\| G\|_{H^1}\le 2 \|  G \|_\cP .
\end{equation}
The B. Davis inequality ( see \cite{sia} ) asserts that 
$$ C^{-1} \| G\|_{H^1}\le  \bE \sup_{ k\le n } |G_k| \le  C \|
G\|_{H^1} . $$
\paragraph{Martingale transforms.}
Let 
$ ( \O , ( \cF_k) , \bP ) $  be a  filtered probability spaces
and   $G = (G_k) $ be an integrable complex valued $ ( \cF_k) $ martingale.
   Define  martingale transforms 
\begin{equation}
T(G) =  \Im \left[\sum_{k=1}^n  w_{k-1} \cdot \Delta G_k
\right]\end{equation}
where $w_k $ is complex valued, $\cF_k $ measurable and $ |w_k | \le 1 . $
Clearly the transform $T$ satisfies  
$$ \| T(G)  \|_{H^1} \le \| G\|_{H^1}   
\quad
{\rm and } 
\quad
\| T(G)  \|_\cP \le  \|  G \|_\cP . $$
This transform, with its unusual imaginary part, 
will be studied in detail in 
Section \ref{mt} where we prove upper $L^1 $ estimates for $T.$
\paragraph{Regular martingales.}
Fix a   filtered probability space
$ ( \O , ( \cF_n) , \bP ) .$ 
We say that $D = ( D_k ) $ is an $\a -$regular martingale  if 
there exist
 $ (\cF _ k) $
adapted 
 sequences  $ (\s_ k)  $ and $ (d_ k)  $ so that 
 \begin{equation}\label{15augper1}\bE _ {k-1} ( \s_k ) = 0 , \quad |\s_k | \le 1 , \quad \bE _ {k-1}
( |\s_k |^2) > \a , \end{equation}
and 
\begin{equation}\label{15augperm2}
 \Delta D_k =  d_{k-1} \s_ k , \quad k \le n . 
\end{equation}
Dyadic martingales are regular, and 
$ \|D \|_{\cP}  \le \|D \|_{H^1} .$

\paragraph{The Hilbert transform.}
The Hilbert transform on $L^2 ( \bT )$ is defined as 
Fourier multiplier by 
$$ H(  e^{in\theta}) = -i {\rm sign} (n)  e^{in\theta} . $$
Let $ h \in H^2_0 ( \bT ) .$ 
Define the complex valued even part of $h $ by 
 $ u( e^{i\theta} ) = \frac12 ( h( e^{i\theta}) + h( e^{-i\theta})) $
The even part of $h $ is sometimes called the cosine series 
of $h. $   The odd part of $ g $ is simply $ v = h - u $.
Since $h $ is analytic and of vanishing mean, 
we recover it from its even part. Indeed
$$  h = u + i H u , $$
where $H$ is the Hilbert transform. 
We also use the following identity 
$$ \int _\bT u( e^{i\theta}) \cos ( \theta ) d \theta 
= - i   \int _\bT v( e^{i\theta} )   \sin ( \theta )  
d \theta . 
$$  
Let  $ f \in L^2 ( \bT ) $ be even, i.e., 
$f(w) =  f( \overline { w} ) , $ then $ g = Hf $ is odd,
i.e., $g(w) = - g (\overline { w} ). $

Let $ h \in H^2_0( \bT )$ and let $ y = \Im h . $
The Hilbert transform recovers $h $ from its imaginary part $y$ ,
we have $ h = -Hy +iy . $ and 
$ \| h \| _2 = \sqrt{2} \| y \|_2 . $ For 
$ w \in \bC ,$ $ | w| = 1 $ we have therefore
$$   \| h \| _2  = \sqrt{2}\| \Im (w\cdot h)\|_2 . $$

\subsubsection*{The Davis and Garsia Inequality for Hardy Martingales}
We review here  results from   \cite{pfxm12}.
The  $L^1 $ norm  of a Hardy martingale $F = (F_k)_{k = 1 }^n $
is  equivalent  to the $L^1 $ norm of its square function. 
This result of J. Bourgain  \cite{b1}, \cite{MR827292}  and a strengthened version
thereof 
 \cite{pfxm12} are  the starting point 
for this work. Specifically we have  
\begin{equation}\label{150812a}
\bE ( \sum_{k = 1 }^n  |\Delta F_k |^2)^{1/2} \le  C \bE  |F|, \end{equation}
and there  are   Hardy martingales  $G = (G_k)_{k = 1 }^n $
and  $B = (B_k)_{k = 1 }^n $ so that 
        \begin{equation}\label{150812b} F = G + B \end{equation} 
\begin{equation}\label{150812c} \bE ( \sum_{k = 1 }^n \bE_{k-1} |\Delta G_k |^2 )^{1/2} 
+ \bE  \sum |\Delta B _k|\le  C \bE  |F|, \end{equation}  
\begin{equation}\label{150812d}
 |\Delta G_k | \le C |F_k | . \end{equation}
We refer to  \eqref{150812a} as 
the square function inequality 
and  to \eqref{150812b} --- \eqref{150812d} as the Davis and Garsia
inequalities for Hardy martingales.
We sketch next a 
unified way of proving 
set of estimates \eqref{150812a} to \eqref{150812d}. See Theorem \ref{hardydavis2}.

\paragraph{Martingale  decomposition. } The joint proof of \eqref{150812a}
to \eqref{150812d} 
is based on the following decomposition obtained in Theorem \ref{9aug123}.  
   To  each  Hardy martingale 
$F = (F_k)_{k = 1 }^n $ there exist  Hardy martingales $G = (G_k)_{k = 1 }^n $
and  $B = (B_k)_{k = 1 }^n $ so that 
\begin{equation}\label{150812g} F = G + B  , \end{equation} 

\begin{equation}\label{150812f} |F_{k-1}| 
+ \frac14 \bE_{k-1}|\Delta B_k|\le \bE_{k-1}  |F_{k}| , \end{equation}   

\begin{equation}\label{150812e} |\Delta G_k | \le C_0| F_{k-1}| . \end{equation}
We  list   consequences  of  \eqref{150812g} --
\eqref{150812e}  leading to  \eqref{150812a} and \eqref{150812c}.
See Theorem \ref{hardydavis2}.
\begin{enumerate}
\item
Taking expectations on both sides of \eqref{150812f} and summing the
resulting telescoping series gives
\begin{equation}\label{150812h}
   \bE  \sum |\Delta B _k|\le 4\bE  |F| .\end{equation}
Subsequently we denote by  $\| B\|_{\cA}$ the norm on the left  hand side
of \eqref{150812h}.
 \item Conditioned to $ \cF_{k-1} $ the martingale differences
$\Delta G_k $ are analytic, square integrable, and of mean zero. Hence
\begin{equation}\label{150812i}
\bE_{k-1} |\Delta G_k |^2  = 2\bE_{k-1} |\Im ( w_{k-1} \cdot \Delta G_k )
|^2 ,\end{equation}
whenever $ w_{k-1}$ is  $ \cF_{k-1} $ measurable and $ | w_{k-1}| = 1
.$
\item
 Theorem \ref{6aug124}
combined with \eqref{150812i} and  \eqref{150812h} gives
\begin{equation}\label{150812j}
 \| G\|_\cP 
\le C  (\| F\|_{H^1} \| F\|_{L^1} )^{1/2}  .
\end{equation}
\item As $F = G + B $ the Burkholder Gundy 
estimate \eqref{150812y}, \eqref{150812j}
  and \eqref{150812h} imply
\begin{equation}\label{150812k}
\|F\|_{H^1}  \le  2 \| G
\|_{\cP} + \| B\|_{\cA} 
\le  C  \| F\|_{L^1}^{1/2} \|F\|_{H^1}^{1/2} . 
\end{equation}
Canceling  the factor $ \|F\|_{H^1}^{1/2} $ gives 
\begin{equation}\label{150812ka}
\|F\|_{H^1} \le C\| F\|_{L^1}. \end{equation}
\item Substitute \eqref{150812ka}  back into  \eqref{150812j} to obtain
\begin{equation}\label{150812l} 
\| G
\|_{\cP} + \| B\|_{\cA} 
\le  C  \| F\|_{L^1} .
\end{equation}
\end{enumerate}
Note that \eqref{150812ka} is the square function estimate 
 for Hardy martingales and that \eqref{150812l} are the Davis and Garsia inequalities.
We view the  Davis and Garsia inequality   as a  lower  bound
for the $L^1$ norm of the Hardy martingale  $F = (F_k)_{k = 1 }^n .$
\subsubsection*{Dyadic  Perturbation and Stability of Davis and Garsia
  Inequalities}
We next discuss the   main results of this paper, 
Theorem \ref{10aug121} and Theorem \ref{hardydavis}.
Fix a Hardy  martingale  $F = (F_k)_{k = 1 }^n $
and  a dyadic martingale  $D = ( D_k )_{k = 1 }^n $ on $ \bT^\bN . $
We say that $F-D $ is a {\it dyadic perturbation} 
of the Hardy martingale  
$F = (F_k)_{k = 1 }^n .$
The  central results   this  paper  are the Davis and Garsia
inequalities for  dyadic perturbations of  
Hardy martingales.

Let  $T$ be the martingale transform 
defined as 
\begin{equation}\label{15082012ee} 
T (H) = \Im \left[ \sum w_{k-1} \Delta H _k  \right],  \quad\quad   w_k =
  (\overline{F_k -D_k})/|F_k- D_k| .
\end{equation}
We obtain in Section \ref{davis}  
the  Davis and Garsia inequalities for 
the  perturbed Hardy martingale $ F- D .$    
That is,  
there exist Hardy martingales $G $ and $B$ so that 
$$ F = G + B , $$ 
\begin{equation}\label{23augb12a}
\|B  \|_\cA \le C \|F-D
\|_{L^1} ,
\end{equation}
and
\begin{equation}\label{17augb40} 
 \|  T (G-D)   \|_\cP  \le C \|F-D
\|_{L^1}^{1/2}\|F -D\|_{H^1}^{1/2} .
\end{equation}
A simple consequence of \eqref{17augb40} is 
\begin{equation}\label{23augb12b}
  \| G  \|_\cP \le C \|F
\|_{L^1} + C \|D\|_{H^1} .  
\end{equation}

See Theorem \ref{10aug121}.

\paragraph{Martingale decomposition.} The proof of \eqref{17augb40}
uses   martingale decomposition  Theorem \ref{hardydavis}:
For each   Hardy martingale  
$F = (F_k)_{k = 1 }^n $ and  dyadic martingale $D = (D_k)_{k = 1 }^n $ 
there exists a Hardy martingale 
$G = (G_k)_{k = 1 }^n $ so that

\begin{equation}\label{15082012a}
 |\Delta G_k | \le C_0| F_{k-1}- D_{k-1}|, \end{equation}

\begin{equation}\label{15082012b} |F_{k-1}- D_{k-1}| 
+ \frac14 \bE_{k-1}|\Delta B_k|\le \bE_{k-1}  |F_{k}- D_{k}| ,\end{equation}  
where
 \begin{equation}\label{15082012c}
B = F - G .\end{equation}
The decomposition  \eqref{15082012a} -- \eqref{15082012c}
implies  the Davis and Garsia inequalities \eqref{17augb40} for dyadic
perturbations
of
Hardy martingales.  See Theorem \ref{10aug121}.   
\begin{enumerate}
\item Take expectations in \eqref{15082012b}  and sum the resulting
telescoping series. This gives
\begin{equation}\label{15082012d}
\|B \|_\cA \le C \|F-D
\|_{L^1}.
\end{equation}
\item Theorem \ref{6aug12}  combined with \eqref{15082012d} and
\eqref{15082012a}
implies that the  martingale transform operator $T$ 
defined in \eqref{15082012ee} 
satisfies 
\begin{equation}\label{15082012e}
\|  T (G-D)   \|_\cP \le C \|F-D
\|_{L^1}^{1/2}\|F -D\|_{H^1}^{1/2} .
\end{equation}
\item By \eqref{15082012e} and the triangle inequality
\begin{equation}
\|  T (G)   \|_\cP \le C \|F\|_{L^1} + C\| D\|_{H^1} 
+ C\| T(D) \|_{\cP} .
\end{equation}
It remains to note that 
$\| G   \|_\cP  = \sqrt 2 \|  T (G)   \|_\cP , $ by analyticity, 
 and that 
$ \| T(D) \|_{\cP} \le \| D\|_{H^1}  $ by regulartity.
\end{enumerate}
\subsubsection*{The Specialization  to  $ D =  \bE_\cD F $}
We  specialize the decomposition \eqref{15082012a}---
\eqref{15082012c} and the estimates \eqref{23augb12a} --\eqref{23augb12b}
to  the case when the dyadic perturbation $D$ is  the conditional expectation 
$\bE _{ \cD} (F) $ of the Hardy martingale $F.$
Then, to each  Hardy martingale  
$F = (F_k)_{k = 1 }^n $ 
there exist  Hardy martingales 
$G = (G_k)_{k = 1 }^n $  and $B = (B_k)_{k = 1 }^n $ so that 
$$ F =  G  + B ,$$  

\begin{equation}\label{15082012g}
\|B \|_\cA \le C \|F-\bE _{ \cD} F
\|_{L^1},
\end{equation}
and 

\begin{equation}\label{15082012h}
\|  T (G- \bE _{ \cD} G )   \|_\cP  \le C \|F- \bE _{ \cD} F
\|_{L^1}^{1/2}\|F \|_{L^1}^{1/2} ,
\end{equation}
where  the martingale transform $T$ operator is given by
\eqref{15082012ee}.
Moreover \eqref{15082012h} gives 
$$ \| G    \|_{H^1} \le 2\| G    \|_\cP \le  C \|F\|_{L^1} .$$
We record  also the  following martingale inequality of independent
interest. See \cite{b1}.
 For each  Hardy martingale $G, $  

\begin{equation}\label{15082012j}
 \| \bE_\cD G \|_{H^1}  \le C\| T ( G - \bE_\cD G) \|_\cP ^{1/4} 
\| G  \|_{\cP}^{3/4}  + C\| G - \bE_\cD G \|^{1/2}_{L^1} \| G \|^{1/2}_{L^1} .
\end{equation}
where $T $ is the martingale transform  that arose in
\eqref{15082012h}. 
See \cite{b1} and the Appendix for 
\eqref{15082012j}.

\subsubsection*{The embedding Theorem revisited}
The embedding theorem  of J. Bourgain \cite{b1}  states that 
$L^1( \bT ) $ is isomorphic to a subspace of  $L^1( \bT )/H_0 ^1( \bT ) . $
The construction of the $L^1$ subspace in  $L^1( \bT )/H_0 ^1( \bT )
$ 
exploited  Hardy martingales. The  $L^1 $ distance 
of a dyadic martingale  to the space of  integrable Hardy 
martingales is the key to  Bourgain 's  embedding \cite{b1}:

There exists 
$\d >0 $ such that for each dyadic martingale $ D = (D_k)_{k = 1 }^n $
on $ \bT ^\bN , $
\begin{equation} \label{16aug121}
\inf \| D - F \|_{ L^1 } > \d   \| D \|_{ L^1 } ,
\end{equation}
where the infimum is taken over all integrable  
Hardy martingales $F .$

The $L ^1$ distance  estimate 
\eqref{16aug121}
results 
from the following inequality of Bourgain  \cite{b1}. (See Section
\ref{davis}.)  
There exists $ C > 0 $ and $ \a > 0 $ so 
that for any Hardy martingale $F = (F_k)_{k = 1 }^n $  
 \begin{equation} \label{16aug122} \| \bE_\cD F \|_{L^1}  \le C \|F-\bE_\cD (F)
\|_{L^1} ^\a\|F
\|_{L^1}^{1-\a}.\end{equation}  
Note that    \eqref{16aug122}  is self improving. It implies  
that there exists $A_0  > 0 $ so that for each 
Hardy martingale,
 \begin{equation}\label{17aug12red1n}
\|  F \|_{L^1}  \le A_0 \|F-\bE_\cD F
\|_{L^1} .\end{equation}
Indeed, for a given $F$ consider separately the cases
$$ \| \bE_\cD F \|_{L^1}
\ge \frac12 \|F\|_{L^1},  
\quad{\rm and } \quad  \| \bE_\cD F \|_{L^1}\le  \frac12 \|F\|_{L^1} $$
In the first case  \eqref{16aug122} gives
\eqref{17aug12red1n} 
by arithmetic. 
In the second case  \eqref{17aug12red1n} follows by
applying triangle inequality 
$$\|  F \|_{L^1}  \le \|F-\bE_\cD F\|_{L^1} +\|\bE_\cD F\|_{L^1} 
 \le \|F-\bE_\cD F\|_{L^1} + \frac12  \|  F \|_{L^1} .          $$
\endproof
\paragraph{Proving \eqref{16aug122}: } Davis and Garsia inequalities 
apply to the proof of \eqref{16aug122}:  
We use  the 
estimates
\eqref{15082012g},  \eqref{15082012h} in combination with 
  \eqref{15082012j}
to  show  that  \eqref{16aug122} holds.
 The martingale transform 
  \eqref{15082012ee}   is the crucial link between the right hand side
  of
\eqref{15082012h} 
  and the left hand side of  \eqref{15082012j}.
\begin{enumerate}
\item Let $F$ be a Hardy martingale  with  $\| F \|_{ L^1 } = 1.$
Assume that 
\begin{equation} \label{16aug12310} 
\|F-\bE_\cD (F)
\|_{L^1}  = \e \end{equation} 
with
 $\e << 1 ,$ 
since otherwise there is nothing to prove.
\item Determine Hardy martingales $ G , B , $ 
satisfying  \eqref{15082012g} and 
\eqref{15082012h} 
and so  that 
 \begin{equation} \label{16aug124}
 F= G + B
 \end{equation}
\item 
By \eqref{16aug12310},  \eqref{16aug124}  and \eqref{15082012g} we have,
 \begin{equation} \label{16aug125}
 \| G \|_{ L^1 } \le  1 +\e, \quad 
 \|G-\bE_\cD G\|_{L^1}  \le  C\e , \quad \| \bE_\cD B \|_{H^1}  \le
 C\e.
 \end{equation}
\item Invoke \eqref{15082012j} and use \eqref{16aug125}
to get
 \begin{equation} \label{16aug126}
\| \bE_\cD F\|_{H^1} \le \| \bE_\cD G \|_{H ^1} + \| \bE_\cD B \|_{H^1}
\le C\| T ( G - \bE_\cD G) \|_\cP ^{1/4} +  C\e^{1/2},
\end{equation}
where $T $ is the martingale transform  operator  arising in  
\eqref{15082012h}.
\item By the  estimate \eqref{15082012h},
 \begin{equation} \label{16aug127}
\| T ( G - \bE_\cD G) \|_\cP \le  C\e^{1/2}.
\end{equation}
\item Combining  \eqref{16aug127}  and \eqref{16aug126}
we get $$\| \bE_\cD F\|_{H^1}  \le C\e^{1/8}, $$
hence $\| \bE_\cD F\|_{L^1}  \le C\e^{1/8}, $ as claimed.
 \end{enumerate}
We derived the interpolatory estimate \eqref{16aug122}
in a straightforward manner 
from three basic estimates:   \eqref{15082012g},
\eqref{15082012h}, 
and  \eqref{15082012j}. 
{\it This } was the  motivation for considering
dyadic perturbations of Hardy martingales 
and their  Davis Garsia decompositions.
\paragraph{Acknowlegement.} It is my pleasure to thank
M.  Schmuckenschl\"ager and P. Wojtaszczyk for many helpful discussion
concerning the topics of  this paper. 

\paragraph{Organization:} 
In Section \ref{mt}  we prepare  general  martingale tools.
In Section \ref{brown}  we prepare  the complex analytic tools.
Section \ref{davis} contains Davis Garsia Inequalities  for dyadic
perturbations of Hardy martingales and  their  applications
to the proof of the embedding theorem.
The Appendix 
contains  the estimates
that relate 
Hardy martingales, cosine
martingales and martingale transforms. 
\section{Martingale Transforms}\label{mt}
Let  
$ ( \O , ( \cF_n) , \bP ) $  be a  filtered probability spaces. 

The following class of martingale transforms plays
 a {\it central}  role in the proof of Bourgain's embedding theorem  \cite{b1}. 
\begin{equation}\label{15augmt1}
T(G) =  \Im \left[\sum_{k=1}^n  w_{k-1} \cdot \Delta G_k
\right]\end{equation}
where $w_k $ is complex valued, adapted, and $ |w_k | \le 1 . $
This section contains  the point-wise  estimates for $T$ 
as needed in Theorem \ref{10aug121} and Theorem \ref{hardydavis}.

\subsubsection*{The Transform Estimate.}
The  spaces  $ L^1, $ $ H^1 ,$ $ \cP $ and  $ \cA $ corresponding to
the  filtration $( \cF_n)$  
are defined in \eqref{15augmt2} and  \eqref{15augmt3}.
\begin{theor}\label{6aug124}
Let   $ F = ( F_k ) $ be a martingale. Let $A > 0 .$  If  the martingale   $ G = ( G_k ) $ 
satisfies 
\begin{equation}\label{24aug1240} |\Delta G_k |\le A |F_{k-1}| ,
\end{equation}
then the  transform $T$ defined as 
\begin{equation}\label{20aug1240x}
T(G) =  \Im \left[\sum_{k=1}^n  w_{k-1} \cdot \Delta G_k  \right]
\quad\quad\text{where}\quad  w_{k-1}=\frac{\overline{F_{k-1}}}{| F_{k-1}|} ,
\quad \quad k \le n \end{equation}
satisfies the point-wise estimate
\begin{equation}\label{20aug1240}
  \| T(G) \|_{\cP}   \le     C  \| F\|_{L^1}^{1/2} 
\|F\|_{H^1}^{1/2}  + C\| B\|_{\cA},
\end{equation}
where $ C = C(A)  $ and
$$ B = F-G .$$ 
\end{theor}
\paragraph{Comments:}
\begin{enumerate} 
\item The hypothesis \eqref{24aug1240}
matches property \eqref{150812e} of the martingale 
decompositions. 
\item  
The  sequence $ ( w_k ) $ defining  $T$ in \eqref{20aug1240x} depends
on the martingale $F$ appearing on the right
hand side of  \eqref{20aug1240}.
\item
The  appearance of the  
$\| F\|_{L^1}^{1/2}$ on the right hand side of 
   \eqref{20aug1240} makes  Theorem \ref{6aug124}  our basic tool.
(Note, $T$ is obviously a contraction on $\cP $ and on $H^1 . $)
\item The unusual imaginary part 
in the definition of  $T $
is the price to  pay for  
having  the  $L^1-$ factor.
\item We will apply   Theorem \ref{6aug124}  only to those decompositions 
$F = G + B $ for which   
$ \| B\|_{\cA} $ is properly under control.
\end{enumerate}
\subsubsection*{Iteration}  
The   following  iteration method \cite{b1} provides the  framework 
 to produce the conditional square function estimates
of  Theorem \ref{6aug124}. 
\begin{theor} \label{iteration}
Let $n \in \bN . $ Given non-negative and integrable 
  $  M_1,  \dots, M_n, $
  $   V_1,  \dots, V_n,$ and  integrable   $  w_1 ,\dots,  w_n$
so that 
the following estimates hold:
\begin{equation}\label{hypo}
\bE ( M^2_{k-1}  +  V_k^2 )^{1/2}  +  \bE w_k 
\le \bE M_k
\quad\text{ for  $1 \le  k \le n .$} \end{equation}
Then 
\begin{equation}\label{sum}
\bE (\sum_{k = 1 }^n V_k ^2 )^{1/2} 
+ 
\bE \sum_{k=1}^n w_k
\le 2 (\bE M_n)^{1/2}(\bE\max_{k\le n} M_k)^{1/2}
\end{equation}
\end{theor}
See \cite{pfxm12} for a proof of   Theorem \ref{iteration}. 
We use it here to get Theorem \ref{6aug124}.
Proposition \ref{8aug121} is needed to establish the hypothesis \eqref{hypo}
of  Theorem \ref{iteration}. 
\begin{prop}\label{8aug121}
Fix a probability space $ (\O , \bP) .$   
Let $A > 0,$  $ z \in \bC ,$ $ w = \overline{z} / |z| . $ 
Assume that 
$$ |g| \le A| z| \quad\text{and}\quad  \int_\O g d\bP =  0 .$$
Then  for  $$y =\Im  (g\cdot  w ) $$ the following estimate holds,
\begin{equation}\label{15augmt4}
  \left( |z|^2 + \a^2 \int_\O  y^2  d\bP \right)^{1/2}  \le \int_\O| z
+g|d\bP,  \end{equation}
where $\a = \a ( A) . $ Consequently, for any  integrable $f$
\begin{equation}\label{10aug1220}  \left( |z|^2 + \a^2 \int_\O  y^2  d\bP \right)^{1/2}  \le \int_\O| z
+f|d\bP + \int_\O|f-g |d\bP .  \end{equation}  
\end{prop}
We first treat the case when $z = 1 $ in a separate Lemma. 
Thereafter we prove the estimates \eqref{15augmt4} and 
 \eqref{10aug1220}.

\begin{lemma} \label{8aug122}
Let $ A > 0 . $ 
Assume that 
$$ |g| \le A  \quad\text{and}\quad  \int_\O g d\bP =  0 .$$
Then $u =\Im (g ) $ satisfies
\begin{equation}\label{15augmt5e}
  \left( 1 + \a^2 \int_\O  u^2   d\bP \right)^{1/2}  \le \int_\O| 1+g|d\bP, \end{equation}
where $\a = \a ( A) . $
\end{lemma}  
\proof  The idea used below is taken from  \cite{b1} p. 695. 
Since $\int g = 0 $,
\begin{equation}\label{15augmt5}
\int_\O |1 + \frac{g}{M}|d\bP \le \int_\O |1 + {g}|d\bP .
\end{equation} 
for any $M \ge 1 . $
Now write $ g = u + i v $ where $u , v $ are bounded real valued. 
Choose  $M = M ( A ) > 2A . $  Then 
rewrite
\begin{equation}\label{15augmt6a} |1+\frac{g}{M}| = 
\left(( 1 + \frac{u}{M})^2  + \frac{v^2}{M^2}\right)^{1/2}.
\end{equation}
Pull out the factor $( 1 + u /M) $ from the term on
the right hand side  of \eqref{15augmt6a}. By arithmetic 
the right hand side  of \eqref{15augmt6a} becomes 

\begin{equation}\label{15augmt6}
\begin{aligned}
 \left( 1 + \frac{u}{M}\right ) \left( 1 + \frac{v^2}{(M+u)^2}
\right)^{1/2}.
\end{aligned}
\end{equation}
Since $ 1 + u / M  > 0 $ and 
\begin{equation}\label{15augmt6b}
\left( 1 + \frac{v^2}{(M+u)^2}\right)^{1/2} \ge 1 + \frac{v^2}{3(M+A)^2}  , \end{equation}
we multiply  $ 1 + u / M $ and the right hand side of 
\eqref{15augmt6b} to  obtain  
\begin{equation}\label{15augmt6c}
|1 + \frac{ g }{M}  | 
 \ge 
1 + \frac{u }{M}   +  \frac{v^2   }{6(M+A)^2}. 
 \end{equation}
Since $ g$ has vanishing mean, we have also 
 $\int  u  = 0 , $ hence taking the expectation in 
\eqref{15augmt6c} gives 
\begin{equation}\label{15augmt6d} 
\int_\O |1 + \frac{g}{M}|d\bP   \ge 1  +  \int_\O \frac{v^2   }{6(M+A)^2}
d \bP. \end{equation}
Finally, since  the right hand side of \eqref{15augmt6d} 
is $ \ge 1 $ we may replace it by its 
square-root and arrive at \eqref{15augmt5e}.

\endproof

\paragraph{Proof  of  Proposition \ref{8aug121}.} 
Scaling and rotation reduces matters the special case of Lemma \ref{8aug122}.
Write $ 1/ z  = ( 1 / |z| )( w) $ where $ w =  \overline{ z} / |z| .$ 
 Then 
$$ 
 \int_\O| z + g|d\bP = |z| \int_\O |  1 + \frac{1}{|z|}  (g\cdot w) | d\bP .  $$
Next put $ y =\Im ( g\cdot w) . $ 
Since $|g|/|z| \le A , $  Lemma \ref{8aug122}   implies 
\begin{equation}\label{15augmt6f} 
 |z| \left( 1 + \frac{\a^2}{|z|^2} \int_\O 
y^2 
d\bP  \right) ^{1/2}  \le \int_\O| z + g|d\bP ,\end{equation} 
where $\a = \a ( A ) . $ 
Note that \eqref{15augmt6f} is just the same as \eqref{15augmt4}.
The triangle inequality gives now \eqref{10aug1220}.
\endproof

\paragraph{Proof of Theorem \ref{6aug124} .} 
Fix $ k \le n . $ Condition on $ \cF_{k-1} $ and put 
$$
z = F_{k-1}, \quad w = \overline{ F_{k-1}}/ |  F_{k-1} | , \quad f = \Delta F_{k} \quad\text{and}\quad g =
\Delta G_{k} . $$
Apply Proposition \ref{8aug121} with the above specification.  Then \eqref{10aug1220} 
implies that
 $$Y_k =\Im ( w_{k-1} \cdot \Delta G_k) , $$ 
satisfies the following estimate
$$( |F_{k-1}|^2 +  \a^2  \bE_{k-1} Y_k^2 )^{1/2}\le 
\bE_{k -1} | F_k |  + \bE_{k -1} | \Delta B_k | . $$
Taking expectation yields 
$$
\bE ( |F_{k-1}|^2 +  \a^2  \bE_{k-1} Y_k^2 )^{1/2}
\le \bE | F_k| +  \bE  | \Delta B_k |   $$
hence Theorem \ref{iteration}    gives
$$\bE ( \sum_{k =1}^n 
\bE_{k-1} Y_k^2 )^{1/2} 
\le  C (\bE | F_n| )^{1/2} \cdot (  \bE \max_{k \le n} | F_k|)^{1/2} +
C   \bE \sum_{k =1}^n | \Delta B_k |   .
$$
 Since 
$T( G ) = \sum Y_{k} $ we get 
$$  \| T(G  ) \|_{\cP}  =  \bE (\sum \bE_{k-1} |Y_{k}|^2  )^{1/2}. $$
The theorem of B. Davis \cite{sia} asserts that  $\bE \max_{k \le n} | F_k| \le C \| F \|_{H^1} .$
This completes the proof.
 \endproof

\subsubsection*{Perturbation} 
We next perturb the  martingales   in 
Theorem \ref{6aug124}  by regular martingales and 
obtain estimates for  the resulting martingale transforms
The perturbations we  consider here are not {\em small} in any sense,
but rather structurally simple. 

Regular martingales on a  filtered probability space
$ ( \O , ( \cF_n) , \bP ) $ are defined in  \eqref{15augper1}, \eqref{15augperm2}.

\begin{theor}\label{6aug12} 
Let $ F = ( F_k )_{k=1}^N $  and 
 $D = (D_k) $ be martingales, so that  
$D = (D_k) $ is  $\a -$regular.
Let 
 $  G = ( G_k )_{k=1}^N $
be a martingale such that
\begin{equation}\label{10aug29}
 |\Delta G_k |\le A |F_{k-1}- D_{k-1}|, \quad \quad k \le N . 
\end{equation} 
Then  the  transform  $T$ given by 
$$ T (H) = \Im \left[ \sum w_{k-1} \Delta H _k  \right],  \quad\quad   w_k =
  (\overline{F_k -D_k})/|F_k- D_k| . $$
satisfies 
$$  \| T(G-D) \|_{\cP}   \le    C  \| F-D\|_{L^1}^{1/2} \|F-
D\|_{H^1}^{1/2} 
+ C \| B \|_{\cA } ,
$$
where $ C = C(A, \a)  $ and 
$$ B = F - G . $$
\end{theor}
Note that the hypothesis \eqref{10aug29} in Theorem \ref{6aug12} 
is matched by the   decompositions 
for Hardy martingales \eqref{15082012a}.

The proof  of Theorem \ref{6aug12} is based on  the iteration
principle  Theorem \ref{iteration}.  
We use Proposition  \ref{8aug125} below to  verify the assumptions
\eqref{hypo}.
Following is an extension of   Proposition \ref{8aug121}.
\begin{prop} \label{8aug125} Fix $ 0<\a \le 1 , C \ge 1 $ and a probability space 
$(\O, \bP) . $  
Let $ \s : \O \to \bC $ satisfy 
$$ \int \s = 0 , \quad |\s |\le 1, \quad \int |\s|^2 > \a . $$
Let $ z \in \bC , w = \overline{z} / |z| , $ and assume   $ g : \O \to \bC $ satisfy
$$ \int g = 0 , \quad |g| \le C |z| .$$ 
Then for any $ b \in \bC ,$ 
$$y = \Im( g - b \s))\cdot w) $$ 
satisfies
$$ \left( |z|^2 + \d^2 \int y^2 \right)^{1/2} \le 
\int | z + g - b\s | , $$
where $ \d = \d( \a , C ) . $ Consequently, for any integrable $f$ 
\begin{equation}\label{11aug1a}
 \left( |z|^2 + \d^2 \int y^2 \right)^{1/2} \le 
\int | z + f - b\s |+ \int | f- g  |     . \end{equation}
  
\end{prop}
\proof Put $ A = 4C / \a . $ Choose $b \in \bC.$ Then we distinguish
between two cases.
\paragraph{Case 1.}
Let $ |b | \le A |z| . $
In that case we have  $$| g - b\s | \le ( C + A ) |z| .$$ 
 Proposition \ref{8aug121} implies  that  $y = \Im( g - b \s)\cdot w $ satisfies
$$ \left( |z|^2 + \d^2 \int y^2 \right)^{1/2} \le \int | z + g - b\s | , $$
with  $ \d =  \d ( \a , C ) $ 
\paragraph{Case 2.}
Let   $ |b | \ge A |z| . $  This case is straightforward since $b$
dominates everything else.
By rotation we assume that $ w = 1 .$
Define the testing function $ m = -\s b / |b| . $ Note that
$$ \int m = 0 , \quad |m |\le 1, \quad \int b \s \overline{m} = -|b| \int |\s|^2  . $$
This gives,
$$
 \int | z + g - b\s | \ge \int ( z + g - b\s ) \overline{m} 
                       \ge |b|\a - \int |g| .
$$
Since $ |g| \le C |z| $ and  $ |b|\ge A |z |  $ 
$$
 \int | z + g - b\s | 
\ge |z| +  |b| A^{-1}\left(A \a
 - C - 1 \right) . $$ 
Recall that we set $ A = 4 C/\a,$ and $ C > 1 . $ Hence 
$A \a  - C - 1  \ge 1 . $ 
 Using again  $ |g| \le C |z| $ and  $ |b|\ge A |z |  $
we have,  that 
$$ 
  |b| \ge \d_0\left(\int | g - b \s |^2
 \right)^{1/2} ,
$$
where $\d_0 =\left( 1 + {C}/A\right)^{-1} . $
Inserting gives, 
$$
  \int | z + g - b\s | \ge       |z |       +    
 \d\left(\int | g - b \s |^2 \right)^{1/2},
$$
where  $ \d =  \d ( \a , C) . $
\endproof
\subsubsection*{Proof of Theorem \ref{6aug12}.}
Fix $k \le N  .$ Condition on $\cF_{k-1}$ and put 
$$
z = F_{k-1} -  D_{k-1} , \quad w = \overline{z}/|z| 
\quad  g = \Delta G_{k} 
\quad\text{and}\quad b \s = \Delta D _k .  $$
Apply Proposition \ref{8aug125} with these parameters. By \eqref{11aug1a} 
   $$ Y_k = \Im (w_{k-1}(\Delta G_{k} - d_{k-1}\s_k))  \quad \text{where} \quad 
 w_{k-1} = \overline{F_{k-1} -  D_{k-1}} /|  F_{k-1} -  D_{k-1}  | 
  $$ 
satisfies
$$( |F_{k-1} - D_{k-1}|^2 + A_0^{-2}\bE_{k-1} |Y_{k}|^2) ^{1/2} 
\le \bE_{k-1}  |F_{k} - D_{k}| + \bE_{k-1}  |\Delta B_{k} | . $$ 
Taking expectations 
gives 
$$\bE( |F_{k-1}- D_{k-1}|^2 + A_0^{-2}\bE_{k-1} |Y_{k}|^2) ^{1/2} 
\le \bE  |F_{k}- D_{k}| + \bE   |\Delta B_{k} |
. $$
By Theorem \ref{iteration}, and the theorem of B. Davis  \cite{sia}
$$\bE (\sum \bE_{k-1} |Y_{k}|^2  )^{1/2}
 \le C  \| F-D\|_{L^1}^{1/2} \|F-
D\|_{H^1}^{1/2} 
+ C \| B \|_{\cA } ,
$$
where $ C = C(A, \a) . $ 
Since 
$T(G - D) = \sum Y_{k} $ we have 
$$  \| T(G-D) \|_{\cP}  =  \bE (\sum \bE_{k-1} |Y_{k}|^2  )^{1/2}. $$

This completes the proof.
\endproof
\section{Brownian Motion and Truncation}\label{brown}
This section contains our complex analytic ingredients
 Our aim is  Theorem \ref{6aug121}. The proofs rely on  
estimates for outer functions  and stopping  time decompositions
for complex Brownian motion.

\subsubsection*{Outer Functions}

Let $H$ denote the Hilbert transform on $L^2 (\bT . ) $ 
Let $ p \in L^\infty ( \bT ) $ be real valued.  Assume that  
$ \log ( 1 - p ) $ is well defined and bounded. 
Then \begin{equation}\label{10aug39}
 q = \exp [\log ( 1 - p ) + i H \log ( 1 - p )] . 
\end{equation}
defines an element in $ H^\infty ( \bT ) .$ 
See Garnett \cite{jg81} for  background.

\begin{lemma} \label{8aug126} 
Let $ p \in L^\infty ( \bT ) $ with  $ 0 \le p  \le 1/2 . $
Then the outer function \eqref{10aug39}
satisfies the following properties
\begin{enumerate}
\item
Then $ p + |q| = 1.$ 
\item Let  $q_1$ be the real part of $q$ and  $q_2$ its imaginary part
  so that $ q = q_1 +i q_2 . $
 Then 
\begin{equation}\label{10aug40a}
\int_\bT q_2 dm = 0,
\end{equation}
and \begin{equation}\label{10aug40b}
 \int_\bT| 1-q_1 | dm \le C_1 \int_\bT p dm , 
\end{equation}
where $C_1 = 7/8 . $ 
\item  If 
$ p $ is even, $p(e^{i\theta} ) =  p(e^{-i\theta} ),$ then 
$q_2 $  is odd,
\begin{equation}\label{10aug40c}
 q_2(e^{i\theta} ) =  -q_2(e^{-i\theta} ). 
\end{equation}
\end{enumerate}
\end{lemma}
\proof
Since 
$$ |q| = \exp [ \log ( 1 - p ) ] = 1 - p , $$
we have  $ p+ |q| = 1  .$ 
Note  $\int H \log ( 1 - p ) = 0 ,$ gives
\eqref{10aug40a}. 
 
We now turn to \eqref{10aug40b}.
Note that 
$$ q_1 = ( 1 - p ) \cos  H(\log (1-p)) , $$
and
$$
1-  q_1 = ( 1 - \cos  H(\log (1-p)) + p  \cos H(\log (1-p)) . $$ 
Since $ 1 - \cos(x) \le x^2 /2 $ we obtain the point-wise estimate
$$ |1-  q_1| \le \frac12 | H(\log (1-p))|^2 +p . $$
By the $L^2 $ estimates for the Hilbert transform,
$$ \int_\bT | H(\log (1-p))|^2 dm \le 2\int_\bT  | \log (1-p))|^2 dm. $$
Invoking  that $0 \le p \le 1/2 $ gives \eqref{10aug40b} as follows 
$$
\int _\bT | \log (1-p))|^2 dm \le  C_1 \int_\bT p^2 dm  \le  C_1 \int_\bT p dm .
$$
If moreover $p $ is even, then by inspection  
$$ q_2 = (1-p) \sin H \log ( 1-p ) $$
 is odd, hence \eqref{10aug40c} holds.
\endproof

\subsubsection*{Brownian Motion}


Let $(B_t)$ denote complex 2D-Brownian motion on Wiener space, 
 and $((\cF_t), \bP ),$ the associated filtered 
probability space. Put
$$ \t = \inf\{ t> 0 : |B_t| > 1 \}. $$
See Durrett \cite{durr}.

The following theorem is our main complex analytic tool.

\begin{theor}\label{6aug121} There exists $C_0 \ge 1 $ so that the
  following holds.
For   $h \in H^1_0( \bT) $  and $ z \in \bC  ,$
let  
\begin{equation} \label{17aug12bm1}
\rho = \inf \{ t< \t : |h(B_t)| > C_0|z| \} ,
 \quad\text{and}\quad 
 g(e^{i\theta}) = \bE ( h(B_\rho) | B_\t = e^{i\theta}) . 
\end{equation}
Then  
$ g  \in H^\infty_0( \bT), $  
\begin{equation} \label{17aug12bm2}
 | g| \le C_0 |z| , \end{equation}
and  for any $b \in \bC $ 
\begin{equation} \label{10aug1249x}
 |z| + \frac14 
\int_{\bT}|h-g|dm
\le \int_{\bT} |z +
h - b \s|dm ,
 \end{equation}
where
$$ \s ( e^{i\theta} ) = { \rm sign} ( \cos( \theta ) ) .$$

\end{theor}
\proof 
By a result of  N. Varopoulos \cite{var1}, 
$g$ defined by \eqref{17aug12bm1} is bounded, analytic with vanishing mean, hence in 
$ H^\infty_0( \bT) .$ See also \cite{jopfxm}.  The upper bound
\eqref{17aug12bm2}
results from \eqref{17aug12bm1}.
We get lower bounds for $\int_{\bT} |z +
h - b \s|dm $ by integrating against testing functions. 
In the   case  $|b| \le 8 |z|  $  we take the   outer functions
of Lemma \ref{8aug126}.  The case    $|b| \ge 8 |z|  $
is straight forward and uses simple exponentials as testing functions. 
\paragraph{Case 1.}
Assume $|b| \le 8 |z| . $ 
Define 
$$ p(e^{i\theta}) = \frac14 \left[ 
\bE ( 1 _A | B_\t = e^{i \theta} )+\bE ( 1 _A | B_\t = e^{-i \theta} )
\right] , \quad\text{where}\quad  A = \{ \rho < \infty \}  .$$
Clearly $p $ is even,  $0 \le  p \le 1/2 ,  $ and
\begin{equation} \label{10aug1250a}
\int_{\bT} p  d m = \frac12 \bP ( A) . 
\end{equation}
By definition of $A = \{\rho < \infty \} ,$
\begin{equation} \label{10aug1250b}
\begin{aligned}
{2} \int_{\bT} |h | p dm & =    \bE  | h( B_\t )  1_A  | \\
                     & \ge C_0 |z| \bP ( A ) .
\end{aligned}
\end{equation}
Using that  $|b| \le 8 |z| , $ in the present case, 
together with  
 \eqref{10aug1250a} and \eqref{10aug1250b} 
we get by the  
triangle inequality 
 \begin{equation} \label{10aug1250c}
\begin{aligned}
\int_{\bT} | z +h -b\s  | p dm &\ge  \frac12 \bE  | h( B_\t )  1_A  | -
(|z|+ |b|)\frac12 \bP(A) \\
& \ge (\frac12 - \frac{9}{2C_0})\bE  | h( B_\t )  1_A  | .
\end{aligned}
\end{equation}
Let    $ q  \in H^\infty ( \bT ) $ be the outer function defined by \eqref{10aug39}.
  Note that   $ q  \in H^\infty ( \bT ) $  is orthogonal to $h .$
By \eqref{10aug40c}   $ q_2  = \Im q $ is an odd function, hence orthogonal to $ \s ,$ and to  constants.
This gives the identities below:
  \begin{equation} \label{10aug1250d}
\begin{aligned}
\int_{\bT} | z +h - b \s | \cdot | q | d m 
&\ge 
\left|\int_{\bT} ( z  + h -b\s ) q
d m \right| \\
&= \left|\int_{\bT} (z+ b\s) q
d m \right|\\
&= \left|\int_{\bT} (z+ b\s) q_1 d m \right|.\\
\end{aligned}
\end{equation}
Recall \eqref{10aug40b}, that
$$ \int_{\bT}|1-  q_1|
d m  \le  C_1 \bP(A)/2 . $$
Hence writing $ q_1 = 1 + (q_1 - 1) $ and  using that  $|b| \le 8 |z| , $ gives 
 \begin{equation} \label{10aug1250e}
\begin{aligned}
\left|\int_{\bT} (z+ b\s) q_1 d m \right| &\ge |z| - \int_{\bT} |(z+
b\s) (1-q_1)| d m \\
& \ge |z| -
 \frac{9C_1}{2} |z| \bP(A) .
\end{aligned}
\end{equation}
Hence  combining   \eqref{10aug1250b}   with  \eqref{10aug1250e} and \eqref{10aug1250d}
gives
 \begin{equation} \label{10aug1250f}
\int_{\bT} | z +h - b \s |\cdot |q| dm \ge |z| - \frac{9C_1}{2C_0} \bE  | h( B_\t )  1_A
| . \end{equation}
Take the sum of  \eqref{10aug1250c} and      \eqref{10aug1250f}.
Since $p + |q| = 1 $ we obtain with 
$ \a_0 =( 1/2 - 9/2C_0 -9C_1/2C_0 ) $ that 
$$\int_{\bT} | z +h -b\s |  dm \ge |z| + \a_0 
\bE  | h( B_\t )  1_A | . 
$$
If $ C_0 > 0 $ is large enough, $ \a _0 > 1/ 4 . $
Since  
$$\int_{\bT}| h-g| dm \le  \bE  | h( B_\t )  1_A |, $$
this gives \eqref{10aug1249x}  in the case $|b| \le 8|z| . $ 

\paragraph{Case 2.}
Next we turn to the case when $|b| > 8|z| . $
This case is  straight forward. The testing functions involved are the simple  exponentials. 
Note first that 
\begin{equation} \label{11aug1f}
 \int_{\bT} | z +h - b \s | dm  \ge \left|\int_{\bT} ( z +h - b \s
  )e^{i\theta} dm \right| = \frac{2|b|}{\pi} . 
\end{equation} 
Next by triangle inequality
\begin{equation} \label{11aug1g}
\int_{\bT} | z +h - b \s | dm  \ge \int_{\bT} |h | dm - (|z|+|b|)
\ge \int_{\bT} |h | dm -\frac{9|b|}{8} .
\end{equation} 
Take a  weighted average of the equation \eqref{11aug1g} and \eqref{11aug1f}, to get
$$ \int_{\bT} | z +h - b \s | dm  
\ge \frac{|b|}{8} +\frac{1}{4}\int_{\bT} |h | dm  .$$
Finally since $|z| \le |b|/8 , $ and $\int_{\bT} |h -g | dm \le
\int_{\bT} |h | dm, $
we get 
$$ \int_{\bT} | z +h - b \s | dm  
\ge |z|  +\frac{1}{4}\int_{\bT} |h -g | dm  .$$

\endproof

\section{Davis and Garsia Inequalities}\label{davis}

 Let $\bT^\bN = \{(x_i)_{i=1}^\infty  \}$ denote 
the countable product of the torus $\bT $
equipped with its product Haar measure.
We return to considering  martingales on  $\bT^\bN$.

\subsubsection*{Davis - Garsia Inequalities for Hardy Martingales revisited}
  We begin, explaining how to get  simultaneously the square function estimate \cite{b1} and the 
Davis-Garsia inequalities \cite{pfxm12} for Hardy martingales.  
We use  the complex analytic 
truncation Theorem \ref{6aug121}.  
and the  transform estimates in  Theorem \ref{6aug124}.  
This proof will be extended further on 
to obtain Davis-Garsia inequalities for perturbed Hardy martingales. 
See  Theorem \ref{hardydavis} below.
\begin{theor}\label{hardydavis2}
For each Hardy martingale  $F = (F_k)_{k = 1 }^n $
\begin{equation}\label{11aug1b} \|F\|_{H^1}  \le  C  \| F\|_{L^1} ,\end{equation}
and there exists  a Hardy Martingale  $G = (G_k)_{k = 1 }^n $
\begin{equation}\label{11aug1c}
 \| G\|_\cP + \|B \|_\cA \le  C \|F \|_{L^1},
\end{equation} 
where
$$ B = F- G $$   
\end{theor}
The    Hardy martingale  $G = (G_k)_{k = 1 }^n $ with the properties 
stated in Theorem \ref{hardydavis2} is obtained in the course of proving the 
 following decomposition.  This construction  carries the complex analytic content of 
the Davies and Garsia inequalities. 
\begin{theor}\label{9aug123} There exists $C_0 > 0 $ so that: 
For any  Hardy martingale $F = (F_k)_{k = 1 }^n $
there is  a  Hardy martingale  $G = (G_k)_{k = 1 }^n $ so that
\begin{equation}\label{11aug10a}
 |\Delta G_k | \le C_0| F_{k-1}| 
\end{equation} 
and 
 \begin{equation}\label{11aug10bb}
|F_{k-1}| 
+ \frac14 \bE_{k-1}|\Delta B_k|\le \bE_{k-1}  |F_{k}| , 
 \end{equation} 
where 
$$ B = F-G .$$
\end{theor}
\proof
We first define the Hardy martingale $ G . $
Fix $k \le n . $ Condition to $\cF_{k-1}  .$
Fix   $x= (x_1, \dots , x_{k-1}) \in \bT^{k-1} $  and $ y \in \bT . $
Put
$$ h(y) = \Delta F_k (x, y )\quad\text{and}\quad
z =  F_{k-1} (x).$$ 
Let 
$$ \rho = \inf \{t < \t : |h( B_t) | > C_0 |z | \} , 
\quad g = \bE( h( B_\rho)| B_\t = e^{i\theta} ) . $$
Then by \cite{var1} $ g \in H^\infty_0 ( \bT) $   and clearly 
 \begin{equation}\label{11aug10b} 
 |g| \le C_0 |z| .\end{equation} 
We apply  Theorem \ref{6aug121} with $ b = 0 $ and  get
 \begin{equation}\label{11aug10c} 
 |z| + \frac14 
\int_{\bT}|h-g|dm
\le \int_{\bT} |z +
h|dm 
. 
\end{equation} 
Put
$$ \Delta G_k(x, y ) = g(y),
\quad\quad\text{and} \quad
\Delta B_k(x , y ) = h(y) - g(y).$$
so that 
$$ \Delta F_k =\Delta G_k + \Delta B_k . $$
Hence \eqref{11aug10b} gives 
$$ |\Delta G_k | \le C_0| F_{k-1}| $$
and by \eqref{11aug10c}
$$ |F_{k-1}| 
+ \frac14 \bE_{k-1}|\Delta B_k|\le \bE_{k-1}  |F_{k}| . $$ 
\endproof
\paragraph{Proof of Theorem  \ref{hardydavis2} .}
Apply Theorem \ref{9aug123} to  $F .$  Let 
$ F = G + B $ be the resulting decomposition into Hardy martingales satisfying
\eqref{11aug10a} and \eqref{11aug10bb}.  Define the rotation
$$ w_{k-1} = \overline{F_{k-1}}/|F_{k-1}| , $$
and the transform
$$  
T(G) =  \Im \left[\sum  w_{k-1} \cdot \Delta G_k  \right] .
$$ 
Since $ |\Delta G_k | \le C_0| F_{k-1}|, $ 
Theorem \ref{6aug124}   implies that  
 \begin{equation}\label{11aug20a} 
 \| T(G) \|_{\cP}   \le     C  \| F\|_{L^1}^{1/2} 
\|F\|_{H^1}^{1/2}  + C\| B\|_{\cA},
\end{equation} 
where $ C = C(C_0) . $
Integrating \eqref{11aug10bb} gives 
$$ \bE |F_{k-1}| 
+ \frac14 \bE |\Delta B_k|\le \bE  |F_{k}| , $$ 
and by
summing the telescoping estimates, one obtains
\begin{equation}\label{11aug20c}
\| B\|_{\cA} \le 4  \|F\|_{L^1} .\end{equation} 
Since  $G$ is a Hardy martingale and  $|w_{k-1}| = 1 $  we have
$$  \bE_{k-1}| \Delta G_k |^2  = 2 \bE_{k-1}|\Im \left[ w_{k-1} \cdot \Delta G_k   \right]  |^2 . $$   
Hence 
\begin{equation}\label{11aug20b} 
  \| G \|_{\cP} = \sqrt{2} \| T(G) \|_{\cP}  .
\end{equation} 
Inserting the estimates  \eqref{11aug20c} and \eqref{11aug20b} 
into equation \eqref{11aug20a}  gives  
 \begin{equation}\label{11aug20d}  
\| G \|_{\cP}   \le     C  \| F\|_{L^1}^{1/2} 
\|F\|_{H^1}^{1/2} .
\end{equation}
It remains to replace in \eqref{11aug20d} the right hand side by $\|F\|_{H^1}.$ 
To this end we use   the Burkholder Gundy inequality
in combination with \eqref{11aug20d}.
$$
\|F\|_{H^1}  \le  2 \| G
\|_{\cP} + \| B\|_{\cA} \\
\le  C  \| F\|_{L^1}^{1/2} \|F\|_{H^1}^{1/2} . 
$$
Cancellation of $\|F\|_{H^1}^{1/2}$  gives the square function estimate,
\eqref{11aug1b} and  with \eqref{11aug20d}  the Davis and Garsia 
inequality
\eqref{11aug1c}  at the same time. 
\endproof
\paragraph{Remarks:}
\begin{enumerate}
\item
The proof of  Theorem \ref{hardydavis2} 
and  Theorem \ref{6aug124} yield general conditions on an 
integrable martingale to be in $H^1 . $
Consider a   martingale 
$F = (F_k)_{k = 1 }^n $  with a decomposition into   $G = (G_k)_{k = 1 }^n $
and  $B = (B_k)_{k = 1 }^n $ so that 
 $F = G + B .$  
Assume that there are $C > 0$  and $ \d > 0 $ so that 
the following conditions are satisfied 

$$ |\Delta G_k | \le C| F_{k-1}| ,$$

$$\bE |F_{k-1}| 
+ \d \bE|\Delta B_k|\le \bE  |F_{k}| ,$$

$$
 \bE_{k-1}|\Delta G_k |^2  \le C \bE_{k-1}|\Im ( w_{k-1} \Delta G_k)  |^2,
\quad\quad w_{k-1}  = \overline{ F_{k-1}}/ |F_{k-1}| .$$
Then 

$$ \|F \|_{H^1 } \le A  \|F \|_{L^1 }, $$
where $ A = A( C , \d )  $
\item
Note also that 

$$\bE |F_{k-1}| 
+ \d \bE|\Delta B_k|\le \bE  |F_{k}| ,$$
and 

$$
 \bE (|F_{k-1}| ^2 + \d \bE_{k-1}|\Delta G_k |^2 )^{1/2} \le  \bE
 |F_{k}|  + C \bE|\Delta B_k|
$$
give 

$$ \|F \|_{H^1 } \le A  \|F \|_{L^1 } .$$ 
\end{enumerate}

\subsubsection*{Dyadic Perturbation and Stability}

The main results of this paper are  
Theorem \ref{10aug121} and Theorem \ref{hardydavis}.  
These theorems  determine to which extent Theorem \ref{hardydavis2} 
is stable under dyadic perturbation. Theorem \ref{10aug122} and its application to 
the embedding theorem \cite{b1} were 
impetus for considering dyadic perturbations of  Hardy
 martingales.

\paragraph{Dyadic martingales.}
We recall the definition of the dyadic  $\s $ algebra on  $ \bT^\bN. $ 
It is  defined by means 
of the independent  Rademacher functions
$$ 
  \s_ k ( x ) = {\rm sign} ( \cos_k( x )) , \quad\quad x = ( x _k ) 
\in  \bT^\bN  , $$
where $\cos_k( x ) = \Re x_k . $ 
Let $\cD $ be the $\s-$algebra on  $\bT^\bN  $ generated by 
$ \s_1,\dots, \s_k, \dots $ .

\paragraph{Stopping time decomposition.} 

We next fix  two martingales,  $F = (F_k)_{k = 1 }^n $ is Hardy and 
 $D = ( D_k )_{k = 1 }^n $ is dyadic. 
Fix $k \le n . $ Condition to $\cF_{k-1}  .$
That is, fix  $(x_1, \dots , x_{k-1})\in \bT^{k-1} .$ 
Put
$$ h(y) = \Delta F_k (x_1, \dots , x_{k-1}, y )\quad\text{and}\quad
z =  F_{k-1} (x_1, \dots , x_{k-1})- D_{k-1} (x_1, \dots , x_{k-1}),$$ 
Let  
\begin{equation}\label{11aug30a} 
 \rho = \inf \{t < \t : |h( B_t) | > C_0 |z | \} , 
\quad g = \bE( h( B_\rho)| B_\t = e^{i\theta} )  ,
\end{equation}
and put
 \begin{equation}\label{11aug30b} 
 \Delta G_k(x_1, \dots , x_{k-1} , y ) = g(y),
\quad\quad\text{and} \quad
\Delta B_k(x_1, \dots , x_{k-1} , y ) = h(y) - g(y), 
\end{equation}
such that  \begin{equation}\label{11aug30c} 
\Delta F_k =\Delta G_k + \Delta B_k . 
\end{equation}
The title of this paper refers to the martingale decomposition defined by  \eqref{11aug30a}  ---  \eqref{11aug30c}. 
We turn now to proving the key properties of the Hardy martingale   $G = (G_k)_{k = 1 }^n $
defined by the stopping times \eqref{11aug30b}.
 
\begin{theor}\label{hardydavis}
For every  Hardy martingale $F = (F_k)_{k = 1 }^n $
and dyadic martingale  $D = ( D_k ) $ the  
 Hardy martingale  $G = (G_k)_{k = 1 }^n $ defined by \eqref{11aug30b}
satisfies the following estimates

 \begin{equation}\label{11aug30d} 
 |\Delta G_k | \le C_0| F_{k-1}- D_{k-1}| 
\end{equation}
and  
\begin{equation}\label{11aug30e} 
 |F_{k-1}- D_{k-1}| 
+ \frac14 \bE_{k-1}|\Delta B_k|\le \bE_{k-1}  |F_{k}- D_{k}|  \end{equation}
where 
$$ B = F - G  .$$ 
\end{theor}
\proof  Fix $k \le n . $ Condition to $\cF_{k-1}  $
and put 
$$ z =  F_{k-1} - D_{k-1} , \quad h =\Delta F_k ,
\quad g =\Delta G_k \quad b \s =  \Delta D_k  .$$ 
By \eqref{11aug30a} $ |g| \le C_0 |z| ,$ hence 
\eqref{11aug30d} holds.
Theorem \ref{6aug121}    gives  
$$ |z| + \frac14 
\int_{\bT}|h-g|dm
\le \int_{\bT} |z +
h -b \s|dm 
. $$
Since 
$$ z +
h -b \s =  F_{k} - D_{k}  \quad \text{and} \quad h-g = \Delta B_k , $$
we  
translate back and get \eqref{11aug30e}. 
\endproof
Following are the consequences of  Theorem \ref{hardydavis}.
For a given  Hardy martingale $F = (F_k)_{k = 1 }^n $
and dyadic martingale  $D = ( D_k ) $ let  $T$ be 
defined as 
\begin{equation}\label{20aug121}
 T (H) = \Im \left[ \sum w_{k-1} \Delta H _k  \right],  \quad\quad   w_k =
  (\overline{F_k -D_k})/|F_k- D_k| . \end{equation}
The next theorem states Davis and Garsia inequalities
for  a perturbed Hardy martingale. In its proof we exploit 
 Theorem \ref{hardydavis}   and   Theorem \ref{6aug12}. 
\begin{theor}\label{10aug121}
For every  Hardy martingale $F = (F_k)_{k = 1 }^n $
and dyadic martingale  $D = ( D_k ) $ the  
 Hardy martingale  $G = (G_k)_{k = 1 }^n $ defined by  \eqref{11aug30b}
satisfies
\begin{equation}\label{11aug40a}
\|B  \|_\cA \le C \|F-D
\|_{L^1} ,
\end{equation}
where
$$  B = F- G , $$
and $T$ defined by \eqref{20aug121} satisfies
\begin{equation}\label{b40}
\|  T (G-D)   \|_\cP \le C \|F-D
\|_{L^1}^{1/2}\|F -D\|_{H^1}^{1/2} ,
\end{equation}
and 
\begin{equation}\label{23augb40}
\| G  \|_\cP \le C \|F\|_{L^1} + C\|D\|_{H^1}. 
\end{equation}
\end{theor}
\proof
Invoke the estimates of  Theorem \ref{hardydavis}.
Taking expectations in  \eqref{11aug30e}. 
gives
$$ \bE |F_{k-1} - D_{k-1} | 
+ \frac14 \bE|\Delta B_k|\le \bE |F_{k} -D_k| . $$ 
Summing the telescoping series gives
\begin{equation}\label{26aug121}
 \bE\sum_{k=1} ^n |\Delta B_k|  \le 4 \bE |F_{n} -D_n|,
\end{equation}
or \eqref{11aug40a}. 
Next use  \eqref{11aug30d} 
and apply  Theorem \ref{6aug12}   to  $T$ 
(defined in \eqref{20aug121})..
This gives
$$  \| T(G-D) \|_{\cP}   \le    C  \| F-D\|_{L^1}^{1/2} \|F-
D\|_{H^1}^{1/2} 
+ C \| B \|_{\cA } \, ,
$$
Invoking \eqref{26aug121} we get 
$$  \| T(G-D) \|_{\cP}   \le    C  \| F-D\|_{L^1}^{1/2} \|F-
D\|_{H^1}^{1/2} ,$$
as claimed. The remaining estimate \eqref{23augb40}
is a simple consequence of the above. 
We get 
first, 
\begin{equation}
\|  T (G)   \|_\cP \le C \|F\|_{L^1} + C\| D\|_{H^1} 
+ C\| T(D) \|_{\cP} .
\end{equation}
Since $ G$ is a Hardy martingale we have 
\begin{equation}
\bE_{k-1} |\Delta G_k |^2  = 2\bE_{k-1} |\Im ( w_{k-1} \cdot \Delta G_k )
|^2 ,\end{equation}
whenever $ w_{k-1}$ is  $ \cF_{k-1} $ measurable and $ | w_{k-1}| = 1
.$ Hence 
$\| G   \|_\cP  = \sqrt 2 \|  T (G)   \|_\cP .  $ 
Note also that for regular martingales
$$ \| T(D) \|_{\cP} \le  \| D\|_{\cP} \le \| D\|_{H^1} . $$ 
This gives  \eqref{23augb40} 
as claimed.
\endproof

\subsubsection*{The special Case $D =  \bE_\cD F$}
We next  {\it specialize}  Theorem \ref{10aug121}.  We fix  a 
 Hardy martingale
$ F  $  and $\bE_\cD F$ its conditional expectation with respect to the 
dyadic  $\s $ algebra.
The resulting martingale transform $T$ is then    
\begin{equation}\label{17augspe1}
 T (H) = \Im \left [ \sum  w_{k-1} \cdot \Delta {H }_k   \right ] ,  \quad\quad   w_k =
  (\overline{F_k  - \bE_\cD  F_k})/|F_k- {\bE_\cD F }_k|
  .\end{equation}
The following theorem records the content of Theorem \ref{10aug121}
in this specialized  setting. 
\begin{theor}\label{10aug122}
For  any Hardy martingale
$ F $ there is a splitting into Hardy martingales 
 $G  $ and $B  $ so that
$ F =G + B , $

\begin{equation}\label{b45aug}
\|B \|_\cA \le C_1 \|F-\bE_\cD F
\|_{L^1} ,\end{equation}
 and  $T$ defined in \eqref{17augspe1} satisfies

\begin{equation}\label{b45}
\|   T (G-\bE_\cD G )  \|_\cP  \le C_1 \|F-\bE_\cD F
\|_{L^1}^{1/2}\|F  \|_{L^1}^{1/2}  .
\end{equation} 
and 

\begin{equation}\label{23augb401}
  \| G  \|_{H^1} \le 2    \| G  \|_\cP \le C \|F\|_{L^1} . 
\end{equation}
\end{theor}
\proof
Apply Theorem \ref{10aug121} to the Hardy martingale
$ F  $  and its conditional expectation $\bE_\cD F .$
Use that $\| \bE_\cD (F) \|_{H^1}  \le C \| F \|_{L^1} . $ 
\endproof

\subsubsection*{An upper  Estimate for $ \bE_\cD G $.}

It remains to complement the Davis Garsia inequalities in 
Theorem \ref{10aug122} with an upper bound  for $ \bE_\cD G $.
For any Hardy martingale 
 $G = (G_k)_{k = 1 }^n $
and    any (!) 
adapted sequence  $W = ( w_k ) $ satisfying   $|w_k | =1 $ the following holds 

\begin{equation}\label{10aug129}
\| \bE_\cD G \|_{H^1}  \le 
 C \| T_W ( G - \bE_\cD G) \|_\cP ^{1/4} 
\| G  \|_{\cP}^{3/4} +C\| G - \bE_\cD G \|^{1/2}_{L^1} \| G \|^{1/2}_{L^1} , \end{equation}
where $T_W $ is the martingale transform  operator
\begin{equation}\label{15aug12mm}
T_W (G - \bE_\cD G ) = \Im \left [ \sum w_{k-1} \Delta_k (G - \bE_\cD G ) \right ] . 
 \end{equation}
See \cite{b1} and the Appendix for  \eqref{10aug129}.
\subsubsection*{The Embedding Theorem revisited}
 J. Bourgain \cite{b1} determines  
a subspace of  $L^1( \bT )/H_0 ^1( \bT )  $
 isomorphic to $L^1( \bT ) .$ 
The construction of such a  subspace relies on the 
following    $L^1-$ distance estimate.

There exists 
$\d >0 $ such that for each dyadic martingale $ D = (D_k)_{k = 1 }^n $
on $ \bT ^\bN , $
\begin{equation} \label{16aug121n}
\inf \| D - F \|_{ L^1 } > \d   \| D \|_{ L^1 } ,
\end{equation}
where the infimum is taken over all integrable  
Hardy martingales $F .$

The embedding theorem will be deduced from the following interpolatory
estimate: 
For any Hardy martingale $F = (F_k)_{k = 1 }^n $  
\begin{equation}\label{10aug128n}
  \| \bE_\cD F \|_{L^1}  \le C \|F-\bE_\cD (F)
\|_{L^1} ^\a\|F
\|_{L^1}^{1-\a}, \end{equation}
 for some $\a >0 . $ See \cite{b1}.
Recall that the estimate   \eqref{10aug128n}  is self improving. 
It implies  that there exists $A_0  > 0 $ so that for each 
Hardy martingale,
 \begin{equation}\label{17aug12red1}
\|  F \|_{L^1}  \le A_0 \|F-\bE_\cD F
\|_{L^1} .\end{equation}
\endproof
\paragraph{Proof that \eqref{17aug12red1}  implies \eqref{16aug121n}.}
The following proof is straightforward, and included  for the  sake 
being definite.
Fix  a dyadic martingale $ D = (D_k)_{k = 1 }^n ,$ 
resolve the inf on the left hand side of \eqref{16aug121n},
thereby  select a Hardy martingale $F_0 $ so that 
$$ \inf \| D - F \|_{ L^1 } \ge \frac12  \| D - F_0 \|_{ L^1 }. $$
If $ \| F_0 \|_{ L^1 }  \le  \| D  \|_{ L^1 }/2 $ we have
$$  \| D - F_0 \|_{ L^1 } \ge \| D  \|_{ L^1 } - \| F_0 \|_{ L^1 } 
\ge \| D  \|_{ L^1 } /2. $$
If conversely  $ \| F_0 \|_{ L^1 }  \ge  \| D  \|_{ L^1 }/2 $ we
proceed by treating separately these two cases: 
\begin{equation}\label{18aug1}
 \| D\|_{ L^1 } \ge 4 A_0  
\|  \bE_\cD F_0 -D\|_{ L^1 }  
\quad
{\rm and}
\quad  \| D\|_{ L^1 } \le 4 A_0  
\|  \bE_\cD F_0  -D\|_{ L^1 } . 
%
\end{equation}
In the first case we write
\begin{equation}\label{18aug2}
 \| D - F_0 \|_{ L^1 } \ge \|  F_0 -\bE_\cD F_0   \|_{ L^1 } - \|
\bE_\cD F_0   - D \|_{ L^1 } . \end{equation}
 Invoke \eqref{17aug12red1} and use the first case in 
\eqref{18aug1}, that is, 
\begin{equation}\label{18aug3}
\|  F_0 -\bE_\cD F_0   \|_{ L^1 }  \ge  (1/A_0) 
\|  F_0   \|_{ L^1 } ,    \quad\quad  
 \| \bE_\cD F_0   - D \|_{ L^1 }  \le  (1 / 4A_0) 
| D   \|_{ L^1 }   
\end{equation}
Since we assumed that   $ \| F_0 \|_{ L^1 }  \ge  \| D  \|_{ L^1 }/2 $ 
we get from \eqref{18aug2} and \eqref{18aug3} that 
 $$\| D - F_0 \|_{ L^1 } \ge  (1 / 4A_0) \| D  \|_{ L^1 } .$$
In the second case of \eqref{18aug1}  write
$$\| D - F_0 \|_{ L^1 } \ge \|\bE_\cD( D - F_0) \|_{ L^1 } 
  = \| D - \bE_\cD F_0 \|_{ L^1 } \ge  (1 / 4A_0)  \| D
\|_{ L^1 } .$$ 

\subsubsection*{Proof of \eqref{10aug128n}. }
The  interpolatory  estimate \eqref{10aug128n}  follows
routinely from   Davis and Garsia inequalities
\eqref{b45} and \eqref{b45aug} of Theorem  \ref{10aug122} 
and  the martingale
estimate \eqref{10aug129}. 
The details of the derivation  are given 
in the following string of remarks. 

Let $F$ be a Hardy martingale and  
apply to it Theorem \ref{10aug122}. 
Let $ F = G + B  $ be the corresponding decomposition. Then
\begin{equation}\label{10aug1212}
\bE_\cD F = \bE_\cD G  +\bE_\cD B . 
\end{equation}

\paragraph{Step 1.} Use \eqref{b45aug} directly to  bound 
$ \bE_\cD B . $ Since $ \Delta ( \bE_\cD B) _k =  \bE_\cD  ( \Delta (  B) _k ) , $ 
we get with \eqref{b45aug},
\begin{equation}\label{10aug1213}
 \| \bE_\cD B \|_{L^1} \le  \| B \|_\cA \le C_1 \|F-\bE_\cD (F)
\|_{L^1} . 
\end{equation}
\paragraph{Step 2.} Use \eqref{10aug129} to bound $ \bE_\cD G . $
 This gives 
\begin{equation}\label{10aug1210}
  \| \bE_\cD G \|_{L^1}  \le
\| T ( G - \bE_\cD G) \|_\cP ^{1/4} 
\| G  \|_{\cP}^{3/4} + C\| G - \bE_\cD G \|^{1/2}_{L^1} \| G \|^{1/2}_{L^1}  , \end{equation}
where the operator $T $ is the one appearing in \eqref{b45}.
\paragraph{Step 3.} 
  Theorem \ref{10aug122} controls the terms on the right 
hand side of \eqref{10aug1210}. Indeed
 \eqref{b45} gives
\begin{equation}\label{10aug1210a}
\|   T (G-\bE_\cD G )  \|_\cP  \le C \|F-\bE_\cD F
\|_{L^1}^{1/2}\|F  \|_{L^1}^{1/2}  . \end{equation}
Moreover by \eqref{23augb401},  
$
\| G  \|_{H^1} \le 2 \| G  \|_{\cP} \le  C\|F  \|_{L^1} .
$
\paragraph{Step 4.} The remaining factor in \eqref{10aug1210}
is $\| G - \bE_\cD G \|_{L^1} .$  By triangle inequality
  $$\| G - \bE_\cD G \|_{L^1}   \le  \| F - \bE_\cD F \|_{L^1}  +  \| B - \bE_\cD B \|_{L^1} . $$
By 
\eqref{10aug1213}, 
 $$
 \| B - \bE_\cD B \|_{L^1} \le 2\| B \|_{\cA}   \le 8  \| F - \bE_\cD F \|_{L^1}  .$$
Hence
\begin{equation}\label{10aug1210c}
\| G - \bE_\cD G \|_{L^1}   \le C \| F - \bE_\cD F \|_{L^1} . \end{equation}
Estimating the  right hand side of \eqref{10aug1210}  using 
 \eqref{10aug1210a}
   and \eqref{10aug1210c},  
gives 
\begin{equation}\label{10aug1214}
  \| \bE_\cD G \|_{L^1}  
\le C \| F - \bE_\cD F \|_{L^1} ^{1/8}   \| F \|_{L^1} ^{7/8} .
\end{equation}
\paragraph{Step 5.}
The identity   \eqref{10aug1212}  and      \eqref{10aug1213},  \eqref{10aug1214}
give 
$$ \| \bE_\cD F \|_{H^1}  \le 
C \|F-\bE_\cD (F)
\|_{L^1} ^\a\|F
\|_{L^1}^{1-\a}, 
\quad {\rm  with} \quad \a  = 1/8. $$ 
 
\endproof

\section{Appendix I. Sine and Cosine Martingales}\label{app1}

In the course of proving  \eqref{10aug128n} we invoked 
\eqref{10aug129}. The proof of \eqref{10aug129} is done in two 
independent propositions concerning estimates for cosine martingales.
Cosine martingales are defined in 
\eqref{14aug12hh}.
Recall that the cosine-martingale  $U =  ( U_k ) $ of 
a Hardy martingale 
$G =( G_k ) $ 
is defined by
\begin{equation}\label{11aug60a}
\Delta U _k (x,y)   = 
\frac12 \left[ \Delta G_k(x,y)+  \Delta G   _k(x,\overline{y})\right], 
\quad    
x \in \bT^{k-1} , \,  y \in \bT . 
 \end{equation}
Proposition \ref{14aug12a} and  Proposition \ref{11aug60b} 
form  the link between  cosine martingale,   $ \bE_\cD G $
and  $ G - \bE_\cD G.$

\begin{prop}\label{14aug12a} Let $G = (G_k) $ be a Hardy martingale, and 
let  $U =  ( U_k ) $ be  its  cosine martingale defined by
\eqref{11aug60a}. Then

\begin{equation}\label{11aug12a}
 \| \bE_\cD G \|_{H^1}  \le C\| U -\bE_\cD U\|^{1/2}_{H^1} \| G\|^{1/2}_{H^1}+ C\| G - \bE_\cD G \|^{1/2}_{L^1} \| G \|^{1/2}_{H^1} 
\end{equation} 
\end{prop}
The   proof of  Proposition \ref{14aug12a}
is in \cite{b1} pp. 700 --702. 
\begin{prop}\label{11aug60b}

    Assume that  $W = ( w_k ) $ is adapted and  $|w_k | =1 ,$
then

\begin{equation}\label{11aug95a}
\|U - \bE_\cD U \|_{\cP} \le 
C \| T_ W ( G - \bE_\cD G) \|_\cP^{1/2} \|G\|_{\cP}^{1/2} ,\end{equation}
where $T_W$ is defined by \eqref{15aug12mm}.
\end{prop}
 The  proof of Proposition \ref{11aug60b}
is in \cite{b1} p. 700.
\subsection*{Randomizing  Martingales}

We give the   proof of  Proposition \ref{14aug12a}
as in \cite{b1} pp. 700 --702. 
%
The
interpolatory estimates \eqref{15aug12xx}
for   the averaging projection $P$
defined in  \eqref{15aug121} are the central 
ingredient.

\paragraph{Randomizing.}

Let $ V = ( V_k) $ be any 
 martingale on $\bT^\bN $. We  define two ways of randomizing  $ V = ( V_k). $
\begin{enumerate}
\item
We associate to   $ V = ( V_k) $ a famly of martingales on 
$\bT^\bN $ parametrized by $ \varepsilon \in \{-1, 1 \} ^{\bN } . $
Put $ v_k = \Delta V _k $ and define

$$ v_k (  x , \varepsilon   ) = 
v_k( x_1^{ \varepsilon  _1} ,\dots,   x_{k}^{ \varepsilon_{k}})
, $$
and 
$$
V( x ,  \varepsilon ) = \sum  v_k (  x , \varepsilon   ) ,
\quad\quad x \in \bT^\bN , \quad \varepsilon \in \{-1, 1 \} ^{\bN } . $$
The partial sums of the series on the right hand side 
form  a familiy of martingales on   $\bT^\bN $ parametrized by 
 $ \varepsilon \in \{-1, 1 \} ^{\bN }  $ so that for each 
fixed  $ \varepsilon , $

\begin{equation} \label{24aug121}
\bE _x | V(x, \varepsilon ) | = \bE _x | V(x ) |  = \| V \|_{L^1} .
%
\end{equation}
and
\begin{equation} \label{24aug122}
 \bE _x(\sum_{k=1}^n|  
v_k (  x , \varepsilon   ) |^2  )^{1/2} = 
\bE _x(\sum_{k=1}^n |  v_k (  x  ) |^2  )^{1/2} = \| V\|_{H^1} .
\end{equation}
\item
Next we associate to  the martingale    $ V = ( V_k) $ 
a family of dyadic martingales parametrized by $ x \in \bT^\bN . $
Put 
$$ d_{k-1} ( x , \varepsilon ) = v_k( 
 x_1^{ \varepsilon  _1} ,\dots,   x_{k-1}^{ \varepsilon_{k-1}} , x_k)
, $$
and form the  familiy of dyadic martingales
$$ D( x ,\varepsilon   )  = \sum   d_{k-1} ( x , \varepsilon ) 
\varepsilon_{k} ,  $$
parametrized by   $ x \in \bT^\bN . $
\end{enumerate} 
\paragraph{Rademacher coefficients.}
The Rademacher coefficients of 
the parametrized dyadic martingales $ D( x ,\varepsilon   )$ 
give rise to averaging projections for $ V = (V_ k ) . $   
\begin{enumerate}
\item
 Given  the dyadic martingales 
$ D( x ,\varepsilon   ) $ the Rademacher coefficienets are defined as 
$$ (R_k ( D ) ( x ) = \bE _{\varepsilon} (  
D( x ,\varepsilon   ) \varepsilon_{k}), $$
where again  $ x \in \bT^\bN  $ is just the parameterindex of  the
family  $ D( x ,\varepsilon   ). $
Remark that $R_k ( D ) = \bE _{\varepsilon} ( d_{k-1} )$
and that   $R_k ( D )$ is determined by the martingale difference 
$ v_k = \Delta V_ k .$ We put 
$$
 P(v_k ) = R_k ( D ) $$
and form the linear extension, 
\begin{equation}\label{15aug121}
 P( V ) = \sum  P(v_k )  = \sum R_k ( D ) . 
\end{equation}
\item Bourgain's version of the Garnett Jones inequality ( see \cite{b1}, 
\cite{durr}, \cite{m} ) 
implies that 
for the Rademacher coefficients 
$$ R_k ( D )  = \bE _{\varepsilon} (  
D(\varepsilon   ) \varepsilon_{k}) $$
of a dyadic martingale 
$$ 
 D( \varepsilon   )  = \sum   d_{k-1} (  \varepsilon ) 
\varepsilon_{k}   $$
there holds 
\begin{equation}\label{2sep121}
\sum  |R_k ( D )|^2 \le C \bE _{\varepsilon} |\sum   d_{k-1} (  \varepsilon ) \varepsilon _k  ) |
\bE _{\varepsilon}  (  \sum |d_{k-1} ( \varepsilon )|^2 ) ^{1/2} 
\end{equation}
\end{enumerate}
\paragraph{Sine martingales.}
Comparing the series representation of $ V( x , \varepsilon) $
and $ D( x ,  \varepsilon   ) $ it is clear that we constructed 
 two different objects,
unless we have  
 further assumptions on the underlying martingale  
$ V = ( V_k) .$  
If $ V = ( V_k) $ is a sine martingale ( see \eqref{14aug12h} ) 
then 
$$    v_k ( x , \varepsilon ) =  d_{k-1} ( x , \varepsilon ) \varepsilon_k $$
and $ D( x ,  \varepsilon   )  = V( x , \varepsilon) .$

\begin{prop} \label{14aug12gj} Let $ V = ( V_k) $ be 
a sine- martingale on  $\bT^\bN . $ 
\begin{equation}\label{15aug12xx}
 \| P(V) \| _{H^1} \le C \| V \|_{L^1}^{1/2} \| V \|_{H^1}^{1/2} .
\end{equation}
\end{prop}
\proof
Since  $ V = ( V_k) $ is  a sine- martingale we have
$    v_k ( x , \varepsilon ) =  d_{k-1} ( x , \varepsilon )
\varepsilon_k $
and 

$$ 
P( v_k)  ( x) = \bE _{\varepsilon}( d_{k-1} ( x , \varepsilon ))
 = \bE _{\varepsilon} ( v_k ( x , \varepsilon )\varepsilon_k ). $$
Hence by \eqref{2sep121},
\begin{equation}\label{2sep122}
\sum  |P( v_k)  ( x)|^2 \le C
\bE _{\varepsilon} | V( x ,\varepsilon   ) 
   |
\bE _{\varepsilon}  (  \sum |v_{k} ( x, \varepsilon )|^2 ) ^{1/2}
.\end{equation}
Applying \eqref{2sep122} shows that  
$$
\| P(V) \| _{H^1} = 
\bE _{x} (\sum  |P( v_k)  ( x)|^2 ) ^{1/2} 
 $$
is bounded by a multiple of
\begin{equation}\label{2sep123} \bE _{x}
\left( \bE _{\varepsilon} | V( x ,\varepsilon   ) 
    |
\bE _{\varepsilon}  (  \sum |v_{k} ( x, \varepsilon )|^2 ) ^{1/2}
\right)^{1/2} .\end{equation}
Next apply the Cauchy Schwarz inequality so that \eqref{2sep123} is 
 bounded by
$$
\left(\bE _{x}   \bE _{\varepsilon}
| V( x ,\varepsilon   ) 
|
\right) ^{1/2}
\left(\bE _{x}   \bE _{\varepsilon}    (  \sum |v_{k} ( x, \varepsilon )|^2
) ^{1/2}
\right)  ^{1/2} .
$$ 
Apply Fubini and 
invoke the identities \eqref{24aug121} and \eqref{24aug121} to get
$$\| P(V) \| _{H^1} \le  C \| V \|_{L^1}^{1/2} \| V \|_{H^1}^{1/2} .
$$   
\endproof

Next we record an application to Hardy martingales. 
Let $G = (G_k) $ be a Hardy martingale 
and  $U =  ( U_k ) $ its cosine martingale 
 given by \eqref{11aug60a}.
\begin{prop} \label{15aug12d}
Let $ Z = ( Z_k ) $ be the  $(\cF_k ) $ martingale 
with difference sequence 
$$ \Delta Z _k =  \bE_{k-1} (\Delta U_k \cos_k )\s_k. $$ 
Then,
\begin{equation}\label{14au1240}
 \|P (Z)  \|_{H^1} \le  C \|P (G-U)  \|_{H^1}, \end{equation}
and consequently 
\begin{equation}\label{14au1240a}
\|P (Z)  \|_{H^1} 
\le C\| G-U  \|^{1/2}_{L^1} \| G-U \|^{1/2}_{H^1}. 
\end{equation}
\end{prop}
\proof Note,  $V =  G-U$ is a sine martingale and  
by the analyticity of $\Delta G_k  , $
\begin{equation}\label{14au1210}
 \bE_{k-1} (\Delta U_k \cos_k )  =   -i\bE_{k-1} (\Delta V_k \sin_k ) .\end{equation} 
Next, \eqref{14au1210} implies
$$|P \Delta Z _k |  =  | P \bE_{k-1} (\Delta V_k \sin_k )  | , $$
and

\begin{equation}\label{14au1210a}
\begin{aligned}   
|P \Delta Z_k | 
& \le     \bE_{k-1} |  P  \Delta V_k | .
\end{aligned}
\end{equation}
Hence \begin{equation}\label{14au1210c}
\|(\sum |P \Delta Z_k | ^2 )^{1/2}  \|_{L^1}\le C  
\|(\sum\bE_{k-1}^2  |  P  \Delta V_k |)^{1/2} \|_{L^1}  .\end{equation}
With  Lepingle inequality \cite{ping}, the right hand side of
\eqref{14au1210c} 
is bounded by $C \|P V  \|_{H^1} , $ that is  
\eqref{14au1240} holds.
Since $V =  G-U$ is a sine-martingale, Proposition \ref{14aug12gj} yields
\eqref{14au1240a}.
\endproof

\subsubsection*{Proof of Proposition \ref{14aug12a} .}

\paragraph{Part 1.}
With the notation of Proposition \ref{15aug12d} we claim that
  \begin{equation}\label{15aug12y} \| \bE_\cD G \|_{H^1}  \le C\| U - \bE_\cD U \|_{H^1} +
\| G-U 
\|^{1/2}_{L^1} \| G- U \|^{1/2}_{H^1}. \end{equation} 
To this end 
let $ Z = ( Z_k ) $ be the 
martingale with difference sequence 
$$ \Delta Z _k = \bE_{k-1} (\Delta U _k \cos_k )\s_k. $$ 
Observe that
 $\bE_\cD G = \bE_\cD U.$
The key identity is 
\begin{equation}\label{15aug12g}
  \bE_\cD G = 
   P(Z)    + (\bE_\cD U  - P(Z) ) , 
\end{equation}
where $P$ is defined in \eqref{15aug121}.
Proposition \ref{15aug12d} 
 readily gives estimates for the first summand  $P(Z)$ in 
\eqref{15aug12g}.
\begin{equation}\label{15aug12x} \| P(Z) \|_{H^1} \le C \| G-U 
\|^{1/2}_{L^1} \| G- U \|^{1/2}_{H^1}. \end{equation} 
Next we turn to estimating $\bE_\cD U  - P(Z) $.  Since  
$\bE_\cD U $ is just even,   we  have  $P( \bE_\cD U) = \bE_\cD U ,$
and
 $$ P(Z) - \bE_\cD U = P( Z - \bE_\cD U) .$$
 By definition $Z$ is a cosine martingale, hence  the operator 
$P$ acts as averaging on $ Z - \bE_\cD U  $ and 
$$ 
   \| P( Z - \bE_\cD U)\|_{H^1}   \le C   \|  Z - \bE_\cD U\|_{H^1} . $$
This, and  invoking
Lepingle inequality  gives
 \begin{equation}\label{15aug12c} \| P(Z) - \bE_\cD U \|_{H^1}  \le C \|  Z - \bE_\cD U\|_{H^1} \le C
\| U - \bE_\cD U \|_{H^1}. \end{equation}
Summing up, \eqref{15aug12y}  
follows from the identity  \eqref{15aug12g} combined with  the estimates   
\eqref{15aug12x}   and  \eqref{15aug12c}. 

%

\paragraph{Part 2.}

By  \eqref{15aug12y}  
\begin{equation}\label{11aug12b}
  \| \bE_\cD G \|_{H^1}  \le
C\| G-U 
\|^{1/2}_{L^1} \| G- U \|^{1/2}_{H^1} + C\| U - \bE_\cD U \|_{H^1} . \end{equation} 
Since 
 $$\bE_\cD G = \bE_\cD U,$$
we have  the  identity
\begin{equation}\label{11aug12c}
G-U =   G- \bE_\cD G  -( U -\bE_\cD U ) , \end{equation} 
which gives immediately
\begin{equation}\label{11aug12d}
 \| G-U \|_{L^1} \le \| G- \bE_\cD G \|_{L^1} + C  \| U -\bE_\cD U \|_{H^1} . \end{equation}  
Next invoke   the (routine)  estimates
 $$ \| G- U \|_{H^1} \le C\| G\|_{H^1} \quad {\rm and } \quad \|U -\bE_\cD U \|_{H^1} \le   C \| G\|_{H^1} .$$ 
Thus, \eqref{11aug12d} and \eqref{11aug12c} imply that
$$\| G-U 
\|^{1/2}_{L^1} \| G- U \|^{1/2}_{H^1} \le 
 C\|  G- \bE_\cD G
\|^{1/2}_{L^1}  \| G \|^{1/2}_{L^1}   + C\| U -\bE_\cD U\|^{1/2}_{H^1} \| G \|^{1/2}_{H^1},  $$ 
hence with \eqref{11aug12b} we obtained      \eqref{11aug12a}. 
\endproof

\subsection*{Estimating Cosine Martingales}\label{app2}
We give the  proof of Proposition \ref{11aug60b}
as in \cite{b1} p. 700.

Let  $L_G^2 (\bT ) $ denote the space of (complex valued) even functions in  $L^2 ( \bT), $
and   $L_U^2 ( \bT ) $ the subspace of  $L^2 ( \bT) $ consisting of
(complex valued) odd functions. 
The space   $ L^2 ( \bT ) $ is the direct sum of the  orthogonal subspaces  
$L_G^2 (\bT ) $ and $L_U^2 ( \bT ) ,$  
$$ L^2 ( \bT ) = L_G^2 (\bT ) \oplus L_U^2 ( \bT ) .$$ 
Recall that 
$ \s ( \theta) = {\rm sign} \cos ( \theta) . $ 
Put 
$ w_0 = 1_\bT , $ $ w_1 = \s, $ 
and choose any    orthonormal system $\{w_k: k \ge 2 \}$ in   
$L_G^2 (\bT ) $ so that  $\{w_k: k \ge 0\}$ is an complete orthonormal
basis for  $L_G^2 (\bT ). $
We next observe that 
in $ L^2 ( \bT )$ an ortho-normal basis is given by the system 
$$\{w_k, H w_k : k \ge 0\}.$$
For the the Hardy space  $ H^2 ( \bT )$  the  analytic system 
$$
\{(w_k +i H w_k) : k \ge 0\}$$
 is an  orthogonal
basis with $\|w_k +i H w_k \|_2 = \sqrt{2} , \, k \ge 1 .$  
\begin{prop}\label{10sep121}
Let   $h \in H_0^2 ( \bT) ,$ and   
$  u(z) 
= (h(z) +h( \overline{z}  ))/2 $ be the  even part of $h.$
Then for $ w,b \in \bC , $ with  $ |w| = 1 , $
$$
 \Im ^2(w\cdot(\la u, \s \ra - b))  + \Re^2 (w\cdot\la u, \s \ra) + 
\int_{\bT} |u - \la u, \s \ra \s |^2 d m 
=\int_{\bT} \Im ^2(w\cdot( h - b \s )) d m
$$
\end{prop}
\proof Fix $ h \in H_0^2 ( \bT) $ and  $ w,b \in \bC , $ with  $ |w| =
1  .$  Clearly by replacing $ h $ by $w h $ and $ b $ by $ w b $ 
it suffices to prove the proposition with  $w = 1 . $
Since $ \int u = 0 $ we have
that $$ u =   \sum_{n = 1}^\infty c_n w_n .$$ 
Apply the Hilbert transform and regroup to get
\begin{equation}\label{10sep122}
 h - b\s= (c_1 - b)\s + ic_1 H\s + \sum_{n = 2}^\infty c_n (w_n +i H
 w_n) . 
\end{equation}
Then, taking imaginary parts gives 
\begin{equation}\label{10sep123} 
\Im (h-b\s) = 
\Im (c_1 - b)\s  + \Re c_1 H \s   + \sum_{n = 2}^\infty  \Im c_n w_n + 
\Re c_n Hw_n . 
\end{equation}
By ortho-gonality the  identity \eqref{10sep123} yields 
\begin{equation}\label{10sep124}
\int_{\bT} \Im ^2( h - b\s ) d m
= \Im ^2(c_1 - b)  + \Re^2 c_1 + \sum_{n = 2}^\infty  |c_n|^2 .
\end{equation}
On the other hand, since  $ \int u = 0 $, $c_1 =\la u, \s \ra, $  and $ w_1 = \s $
we get 
\begin{equation}\label{10sep125}
\int_{\bT} |u - \la u, \s \ra \s |^2 d m
 =  \sum_{n = 2}^\infty |c_n|^2 .
\end{equation} 
Comparing the equations \eqref{10sep124}  and \eqref{10sep125}
completes the proof. 
\endproof


We use below  some arithmetic, that we isolate first.
\begin{lemma}\label{10sep126}
Let  $ \mu , b \in \bC $ and
 \begin{equation}\label{11aug65d}
 |\mu| +\dfrac{|\mu-b|^2}{|\mu|+|b|} = a. 
\end{equation}
Then for any  $w \in \bT ,$ 

\begin{equation}\label{11aug60e}
( a - |b|)^2 \le 4( \Im^2 (w\cdot(\mu-b))
+ \Re^2 (w\cdot \mu) )         .
\end{equation}
and 
\begin{equation}\label{11aug60d}
|\mu-b |^2 \le 2(a^2 - |\mu|^2). 
\end{equation}

\end{lemma}
\proof 
By rotation invariance it suffices to prove 
\eqref{11aug60e} for $w = 1 .$ 
Let $ \mu = m_1 + i m_2 $ and $ b = b_1 + i b_2 . $
By definition \eqref{11aug65d}, we have 
$$ 
a - |b| = \dfrac{ |\mu| ^2 - |b|^2 + | \mu - b |^2}{| \mu| + |b |}.
$$ 
Expand and regroup the numerator
\begin{equation}\label{17sep121}
|\mu| ^2 - |b|^2 + | \mu - b |^2 
 = 2m_1 ( m_1  - b_1 ) + 2 m_2 ( m_2 - b_2 ) .
\end{equation}
By the Cauchy Schwarz inequality, the last term in \eqref{17sep121}
is bounded by 
$$ 
2 ( m_1^2 + ( m_2 - b_2 )^2 )^{1/2}(  m_2^2 + ( m_1 -b_1)^2)^{1/2}  
$$
Note that  $m_1 = \Re \mu $   and  $m_2 - b_2 = \Im ( \mu - b ) . $ 
It remains to observe that 
$$
(  m_2^2 + ( m_1 -b_1)^2)^{1/2}  \le |\mu | + |b| .$$
or equivalently 
$$
 m_1^2+ m_2^2  - 2 m _1 b_2 + b_1^2  \le  |\mu |^2 +2   |\mu| |b|+
 |b|^2, $$
which is obviously true. 

Next we turn to verifying \eqref{11aug60d}.  We have
$a^2 -|\mu|^2 = (a+|\mu|)(a -|\mu|)$ hence
\begin{equation}\label{17sep124}
a^2 -|\mu|^2 = 
\left[ 2|\mu| +\dfrac{|\mu-b|^2}{|\mu|+|b|}
\right]\dfrac{|\mu-b|^2}{|\mu|+|b|} .\end{equation}
In view of \eqref{17sep124} we get   \eqref{11aug60d} by showing that  
\begin{equation}\label{17sep123}
 2|\mu|^2 +2|\mu||b| +  |\mu-b|^2 \ge \frac12 ( |\mu|+|b| )^2 .
\end{equation}
The left hand side of \eqref{17sep123} is larger than 
$|\mu|^2+|b|^2$ while the right hand side of   \eqref{17sep123} is
smaller $|\mu|^2+|b|^2 .$ 
\endproof

We  merge the inequalities of Lemma \ref{10sep126} 
with the identiy in Proposition \ref{10sep121}


\begin{prop}\label{25a121} There exists $C_0 > 0 $ so that the
  following holds.
 Let $ w , b \in \bC , $  with $|w| = 1 ,$ $h \in H_0^2 ( \bT) ,$
 let $u $ be the even part of $h $ and put 
 $$  | \la u, \s \ra| +\frac{| \la u, \s \ra-b|^2}{|b|+| \la u,
  \s \ra|}   = a .$$
Then 

\begin{equation}\label{22aug12b}
\int_{\bT} |u - b \s |^2 dm(y)
\le   8(a^2 -|\la u, \s \ra|^2)+ \int_{\bT} |u - \la u, \s \ra \s |^2
dm(y).
\end{equation}
and

\begin{equation}\label{22aug12}
( a - |b| )^2 + \int_{\bT} |u - \la u, \s \ra \s |^2
dm(y) \le 8 
\int_{\bT}  \Im ^2(w\cdot( h - b \s )) dm(y) .                           
\end{equation}  

\end{prop}
\proof
Put 
\begin{equation}\label{25a123}
 J^2 = \int_{\bT}  \Im ^2(w\cdot( h - b \s )) dm(y) .  
\end{equation}
The proof  explpoits the  basic identities for the integral $J^2 $ and 
$\int_{\bT} |u - b \s |^2 dm(y) $ and intertwines them   with  the
arithmetic    \eqref{11aug65d}  -- \eqref{11aug60e}.
\paragraph{Step  1.} Use  the  straight forward identity,
\begin{equation}\label{22aug12a}
\int_{\bT} |u - b \s |^2 dm(y)
= |\la u, \s \ra   -b |^2 + \int_{\bT} |u - \la u, \s \ra \s |^2
dm(y).
\end{equation}
Apply \eqref{11aug60d},  so that  
$$  |\la u, \s \ra   -b |^2 \le 8(a^2 -|\la u, \s \ra|^2 ) ,$$
hence by \eqref{22aug12a} we get \eqref{22aug12b},
\begin{equation*}
\int_{\bT} |u - b \s |^2 dm(y)
\le   8(a^2 -|\la u, \s \ra|^2)+ \int_{\bT} |u - \la u, \s \ra \s |^2
dm(y).
\end{equation*}
      
\paragraph{Step  2.} 
Proposition \ref{10sep121} asserts that 
\begin{equation}\label{22aug12c}
 \Im ^2(w\cdot(\la u, \s \ra - b))  + \Re^2 (w\cdot\la u, \s \ra) + 
\int_{\bT} |u - \la u, \s \ra \s |^2 dm(y) = J^2 .
\end{equation}
Apply \eqref{11aug60e}  with $ \mu = \la u, \s \ra  $
to the left hand side in \eqref{22aug12c}, and
get \eqref{22aug12},
\begin{equation*}
( a - |b| )^2 + \int_{\bT} |u - \la u, \s \ra \s |^2
dm(y) \le 8 J^2 .
\end{equation*}  

\endproof

\subsubsection*{Proof of Proposition \ref{11aug60b}. }

 Let $ \{g_k\} $  be the martingale difference sequence of 
the Hardy martingale $G = (G_k) ,$ and 
let    
 $ \{u_k\} $  be the martingale difference sequence of
the  associated cosine martingale  
$U = (U_k) .$ 
Note  $$   \bE_{\cD} ( u_k ) =  \bE_{\cD} \bE_{k-1}( u_k  \s_k ) \s_k
.$$

\paragraph{Step 1.} 
Put
$ b_k  = \bE_{\cD} \bE_{k-1}( u_k  \s_k ),$ and put
$$
Y^2 = \sum _{k = 1 }^\infty | \bE_{k-1} (u_k  \s_k) |^2 \quad\text{and}\quad
Z^2 = \sum_{k = 1 }^\infty  |b_k| ^2 . $$
We have 
$$
 \bE \bE_\cD (  \sum_{k = 1}^\infty |\bE_{k-1}( u_k  \s_k )|^2
)^{1/2}  \ge \bE  (  \sum_{k = 1}^\infty |\bE_\cD\bE_{k-1}( u_k  \s_k )|^2
)^{1/2} .
$$
hence 
\begin{equation}\label{28aug70j}
\bE ( Y ) \ge \bE ( Z) .
\end{equation} 

\paragraph{Step 2.} 
Since  
 $
\bE_{\cD}(g_k) = \bE_{\cD}(u_k) , $
the square of the conditioned square functions of 
$ T_W (   G -  \bE_{\cD} G ) $  is
\begin{equation}\label{28a12x} 
\sum  \bE_{k-1} |\Im (w_{k-1}\cdot( g_k  - b_k \s_k ))|^2 .
\end{equation}
\paragraph{Step 3.} 
The martingale differences of 
$ U - 
\bE_{\cD}( U)$  is  $ \{ u_k -b_k \s_k\} .$ 
The square of its  conditioned square functions of 
\begin{equation}\label{28a1245}
\sum \bE_{k-1}| u_k - b_k \s_k |^2 .
\end{equation}
Following the pattern of \eqref{11aug65d} define 

$$ a_k  = |\bE_{k-1}( u_k \s_k )| +\frac{|\bE_{k-1}( u_k  \s_k )   -b_k|^2}{ 
|\bE_{k-1}( u_k \s_k )| +   |b_k|} 
 , $$
and   
$$v_k =u_k - \bE_{k-1}( u_k  \s_k )\s_k , \quad\quad  r_k ^2 = \bE_{k-1}| v_k  |^2 .$$ 
 By  \eqref{22aug12b} 

\begin{equation}\label{11aug70a}
\begin{aligned}
&\bE_{k-1}| u_k - b_k \s_k |^2 
\le  8 ( a_k^2  +  r_k ^2  - |\bE_{k-1}^2( u_k  \s_k ) |     ).
\end{aligned}
\end{equation}

\paragraph{Step 4.}
Define
$$ X^2 = \sum_{k = 1 }^\infty a_k ^2 + r_k^2 ,$$
then $ X \ge Y $ and
 \begin{equation}\label{28aug70g}
 \| U - \bE_\cD (U) \|_\cP \le \sqrt{8}  \bE (X^2 -Y^2)^{1/2} \le C( \bE( X - Y ))^{1/2} ( \bE( X + Y ))^{1/2} . 
\end{equation}
 The second factor $\bE( X + Y $ 
in \eqref{28aug70g} is simply   bounded as 
\begin{equation}\label{28aug70h}
 \bE (X+Y) \le C \|U  \|_{\cP} \le  C \|G  \|_{\cP} .
\end{equation}
\paragraph{Step 5.}
Next we turn to  estimates for  $\bE( X - Y ).$ 
First recall 
$\bE( X - Y ) 
\le 
 \bE( X - Z ) ,$
and by triangle inequality
$$
X - Z \le (\sum_{k = 1 }^\infty (a_k - |b_k|)^2 + r_k^2)^{1/2} 
.$$
By \eqref{22aug12} 
\begin{equation}\label{11aug70b}
\begin{aligned}
&( a_k - |b_k|)^2 +  r_k^2 
\le 8 \bE_{k-1} |\Im (w_{k-1}\cdot( g_k  - b_k \s_k ))|^2 .
\end{aligned}
\end{equation}
Combining
\eqref{28aug70g}
---
\eqref{11aug70b}
 gives
\begin{equation}\label{11aug70k}
\bE( X - Y ) 
\le   \bE( X - Z ) 
 \le C \| T_W(G  - \bE_\cD G )\|_\cP .
\end{equation}
\endproof

\section{Appendix II.  The Transfer Operator}\label{app3}
We transfer the the $L^1$ distance estimate from the infinite 
torus to $\bT .$  See \cite{b1} p. 697.   

The space of integrable Hardy amrtingale is denoted by 
 $ H^1 (\bT^\bN ).$
We showed  that there exisits $ A_0  > 0 $ so that 
\begin{equation}\label{13sep125}
 \| D \| _{L^1 } \le A_0   \| F -  D \| _{L^1 } , \end{equation}
whenever $F \in  H^1 (\bT^\bN )$  and $D$ is a dyadic martingale on
 $ \bT^\bN . $
We construct a diffuse  sigma-algebra $ \Sigma $ on $\bT  , $
so that 
$$ \|h \| _{L^1( \bT )  } \le A \| f-h \|  _{L^1( \bT )  } , $$
whenever $ h \in L^1 ( \Sigma) $ and $ f \in H^1 ( \bT ) , $
or equivalently, 
$ \|h \| _{L^1 ( \bT ) } \le A \| h \|  _{L^1( \bT ) / H^1_ 0 ( \bT )
}$
for $h \in L^1 ( \Sigma) .$   
J. Bourgain  \cite{b1} p. 697 determined a bounded linear operator  $J:  L^1  ( \bT )  \to 
L^1 ( \bT^\bN ) $ so that the restrictions to $  L^1 ( \Sigma) $
respectively to $  H^1 ( \bT )   $  satisfy the following conditions.
\begin{enumerate}
\item  The restriction of $J $  to $  L^1 ( \Sigma) $
is an embedding, There exists $A_1 > 0 $ so that
\begin{equation}\label{13sep121} 
 \frac{1}{A_1} \|h \| _{L^1 ( \bT ) }  \le   \| J h \| _{L^1 } \le  A_1  \|h \|
_{L^1( \bT )  }  , \quad\quad h \in  L^1 ( \Sigma) .
\end{equation}
\item The conditional 
expectation operator $ \bE_\cD $ onto the
  subspace of dyadic martingales in 
$L^1 ( \bT^\bN ) $  acts as a
  small perturbation of the identity on $ J (  L^1 ( \Sigma)), $
\begin{equation}\label{13sep122}   \| H -  \bE_\cD H  \|_{L^1 } \le  \frac{1}{4 A_0 } \| \bE_\cD H
\|_{L^1 } \quad\quad H \in  J (  L^1 ( \Sigma)) . 
\end{equation}
\item The restriction of $J $ to $  H^1 ( \bT )   $ maps into the
  space of integrable Hardy martingales, 
 \begin{equation}\label{13sep123}     
J (  H^1 ( \bT )) \sbe  H^1 ( \bT^\bN  ).
\end{equation}
\end{enumerate}

We use the bounded operator $J$   
satisfying \eqref{13sep121}---\eqref{13sep123}
to prove that    $L^1(\Sigma) $ embeds as a closed linear subspace
into  $L^1(\bT) /H^1_0(\bT) . $
\begin{theor}\label{13sep124}
There exists 
$A > 0 $ so that for each $h \in L^1(\Sigma)$
$$
 \| h\|_{L^1 (\bT) } \le A \| h\|_{L^1 (\bT) / H^1
  (\bT)}  
. $$
\end{theor}
\proof 
Let $ h \in   L^1(\Sigma)$ and $ f \in  H^1
  (\bT) . $ Put $H = Jh , $ $ F = Jf $ and $ D = \bE _\cD H . $
By  \eqref{13sep123},  $ F $ is a Hardy 
martingale. Since   $D$ is dyadic by definition, 
\eqref{13sep125} gives 

\begin{equation}\label{15sep1210}  
\| D \| _{L^1 } \le A_0  \|F -  D \| _{L^1 } . \end{equation} 
Write $ F- D = (F - H) + (H - D) .$ 
Since $ D =  \bE _\cD H$  we get  from \eqref{13sep122} 
and \eqref{15sep1210} that  
\begin{equation}\label{15sep121} 
 \| D \| _{L^1 } \le A_0  \|F -  H \| _{L^1 }   +   \frac{1}{4  }\| D
 \| _{L^1 }. 
\end{equation} 
Since 
$F -  H = J ( f - h ) , $  
 we get 
\begin{equation}\label{15sep122}\| D \| _{L^1 } \le 2 A_0   \|J \| \cdot \|f -  h \| _{L^1 (\bT)}  . \end{equation}
Next use that $J $ is an embedding    
of $ L^1 ( \Sigma ). $ 
By \eqref{13sep121} and \eqref{13sep122} we have  
 \begin{equation}\label{15sep123}\| h \|_{L^1 ( \bT )  } 
\le 
2 A_1  \| D  \|_{L^1
} , \end{equation}
since   $ h \in   L^1 ( \Sigma ) $
and  $ D =  \bE_\cD Jh. $ 
Combining \eqref{15sep122} and \eqref{15sep123} gives
$$ \| h \|_{L^1 ( \bT )  } \le  A  \| f - h \|_ {L^1 ( \bT )  } . $$
where $ A = 4 A_1 A_0 \|J\| . $
\endproof

\subsubsection*{The space $L^1(\Sigma ) $ .}
The Fejer kernels $ F_a ,$  $ a \in \bN $  on $ \bT $  
are defined as 
$$ F_a  (z) = \sum _{|j | \le a} ( 1 - \frac{|j|}{a + 1 } )  z^{ j} ,
\quad\quad z \in \bT . $$
We let  $ \s ( z) = {\rm sign} (\Re z) .  $
Define inductively the sigma algebra $\Sigma$ on $\bT . $

\paragraph{Step 1.} Let  $A_0 > 0 $ be the constant appearing in \eqref{13sep125}.
 
Fix $\e > 0  $  where $ \e = \e ( A_0 ) \le (100 A_0)^{-1} . $ Let $ n_1 = 1 . $  .  Select $ a_1 \in \bN $
so that $$ \| s_1  - \s \|_1 < \e/2,  $$ 
where   $ s_1 = K_{a_1 } *\s  .$  Put 
$$ E_1 = \{ k_1 \in ( -a_1 , a_1 ) \cap \bZ \} .$$
  Having defined integers $ n_1 < \dots  < n_{m } $ and
 $ a_1 < \dots  < a_{m } .$
Form 
$$E_m = \{ \sum_{i = 1 } ^m k_i n_i : k_i \in (- a_i , a_i ) \cap \bZ  \}.   
$$
\paragraph{Step m +1 .}
Choose $ n_{m+1} \ge n _m $ so that
\begin{equation}\label{14sep122}
4^m |j| \le    n_{m+1}   ,\quad\quad     j \in E_m  .
\end{equation}
Select  $ a_{m+1} \ge a _m $ so that 
$$\| s_{m+1} - \s \|_1 \le 2^{-m}\e . $$
where $ s_{m+1} =   K_{a_{m+1 }} *\s  $ and  $ \s ( z) = {\rm sign} (\Re z) .  $
\paragraph{The conditional expectation $\bE_\Sigma  $.}
Define $\Sigma$ to be the $\s-$algebra on $\bT $ generated by the
sequence of Rademacher functions
$$  \s (z^{n_k}  ) , \quad z \in \bT . $$ 
Define the non-negative kernel as 
$$  B( z, \zeta ) = \prod_{k = 1 }^\infty ( 1 +  \s( z^{ n_k }
) \s ( \zeta ^{  n_k } )) , \quad \quad z, \zeta \in \bT .
$$ Let $\bE_\Sigma  $ be the conditional expectation operator acting on $L^1 (\bT)$ onto  
$L^1(\Sigma). $ It  is an  integral operator
with kernel  $B( z, \zeta ) , $
\begin{equation}\label{14sep123}
\bE_\Sigma (g) ( z ) = \int_\bT B ( z, \zeta )  g( \zeta ) dm ( \zeta
), \quad g \in L^1 (\bT) . \end{equation}
\paragraph{Products of Fejer kernels}
Put 
$ E = \bigcup E_m. $ 
Let $ (a_k) $ and   $(n_k)$ be the sequence given in the construction
of $\Sigma .$ 
Form the pointwise products of  Fejer kernels 
$$
K( \zeta ) = \prod_{k = 1 }^\infty    F_{a_k} ( \zeta ^{ n_k } ) . $$
We have the following relations for the Fourier coefficients of $K ,$ 

$$ \hat K (0) = 1 ,\quad \text{and} \quad\{m \in \bZ :  \hat K (m) \ne 0 \}
\sbe E , $$
$$
\hat K ( \sum_{k=1}^\infty b_k n_k ) = \prod_{k = 1 }^\infty \hat F_{a_k} ( b_k
). $$
Hence  the Fourier expansion of   $ K $ is as follows,
$$ 
 K( \zeta  )  =  
\sum_{k = 1 }^\infty  \sum_{b_k = -a_k  ,\, b_k \ne 0  }^{a_k} \cdots \sum_{b_1 = -a_1}^{a_1}
\prod_{j = 1 }^k \hat F_{a_j} (b_j) \zeta ^{n_j b_j  } .
$$
 Let $L^1_E = \{ f \in L^1 ( \bT) : \hat f (m) =0 \text{ for } m
\notin E\} . $ We observed  that $ K \in   L^1_E .$
\paragraph{Embedding $L^1(\Sigma)$  into  $L^1 _E .$}
The next proposition identifies the integral kernel of the conditional 
expectation operator $\bE_\Sigma $ after its  convolution with $K.$

\begin{prop} Let  
$  Rg  = K * (\bE_\Sigma g) ,$ $ g \in L^1 ( \bT ) ,$
and put 
$$
A(  z , \zeta  ) = \prod_{k = 1 }^\infty 
( 1 +  s_k(z^{ n_k }) \s ( \zeta^{  n_k } )), \quad\quad z, \zeta \in
\bT .
$$ 
Then
\begin{equation}\label{14sep124}
 Rg ( z ) = \int_\bT A ( z, \zeta )  g( \zeta ) dm ( \zeta ),
\quad\quad z \in \bT . 
\end{equation}
\end{prop}
\proof
We show   that 
$$ K {*_z} B( z, \zeta )   = A( z, \zeta) , $$
where the convolution is taken with respect to the $ z$  variable. 
To this end we observe that for fixed $ \zeta \in \bT $ 
the following identities hold.
$$
\begin{aligned} 
 K {*_z}B (  z, \zeta  )  
&=
 \sum_{k = 1 }^\infty  \sum_{b_k = -a_k  ,\, b_k \ne 0  }^{a_k} \cdots \sum_{b_1 = -a_1}^{a_1}
\prod_{j = 1 }^k   \widehat {F_{a_j}*\s  (b_j)} (z ^{ n_j b_j } ) \s (
  \zeta ^{  n_j } ) \\
&=
 \sum_{k = 1 }^\infty  \sum_{b_k = -a_k  ,\, b_k \ne 0  }^{a_k} \cdots \sum_{b_1 = -a_1}^{a_1}
\prod_{j = 1 }^k  \hat s_j  (b_j) ( z^{ n_j b_j } ) \s ( \zeta ^{ 
  n_j } ) \\
&=
 \prod_{k = 1 }^\infty ( 1 +  s_k(z^{ n_k }) \s ( \zeta^{  n_k } )) .
\end{aligned}
$$
\endproof

\begin{prop}\label{small1}
On $ L^1 ( \Sigma ) $ convolution by $K$ is a small perturbation of 
the identity, 
$$ \|K * h - h \|_{L^1 (\bT )} \le \e \| h \|_{L^1  (\bT )}  \quad\quad h \in  L^1 ( \Sigma ). $$
The operator $  Rg ( z ) = K * (\bE_\Sigma g) $ satisfies 
$$ \| R g - \bE_\Sigma g \|_{L^1  (\bT )}  \le  \e  \| g \|_{L^1 (\bT ) }  , \quad
\quad g \in L^1(\bT) . $$
\end{prop}
\proof In view of the integral representations
 \eqref{14sep123} and \eqref{14sep124}
it suffices to prove that 

$$\sup_{\zeta \in \bT } \int_\bT  | A(z,\zeta)  - B(z,\zeta)  | d m(
z ) \le \e .$$
To this end fix $ \zeta \in \bT ,$ put 
$ \tau_k = \s ( \zeta ^{n_k} )
. $  Let  $ A_0 = B_0 = 1 ,$ 
and for  $ j \in \bN $ put 
$$A_j ( z ) = \prod_{k = 1 }^j ( 1 +  s_k(z^{ n_k }) \t_k) , 
\quad\quad B_j ( z ) = \prod_{k = 1 }^j ( 1 +  \s(z^{ n_k }) \t_k).$$
Rewrite the difference  $A_j ( z )- B_j ( z )$ 
  of the kernels as follows

$$
 A_{j-1}  ( z ) \t_j (s_j( z^{ n_j} ) - \s ( z^{ n_j}  )) 
+     (A_{j-1}  ( z)   - B_{j-1}  ( z )) 
\t_j ( 1 + \s( z^{ n_j})) .
$$
Next take absolute values and exploit that $ ( z^{n_k} ) $ is an
almost independent sequence. Since $   A_{j-1}  \ge 0 ,$ 
$$ \int_\bT A_{j-1}  ( z ) dm(z)  = 1 \quad\text{ and }\quad 
\int_\bT(1 \pm \s( z )) dm(z)  = 1 ,$$
and invoking \eqref{14sep122} it is easy to see that 
$$
 \int_{ \bT} |A_j ( z )- B_j ( z )| dm(z) $$
is bounded by 
$$ 
( 1 + \d_j) \int_{ \bT} |A_{j-1} ( z ) - B_{j-1} ( z ) | dm(z)  + (
1 + \d_j) 
\int_\bT | s_j( z) - \s ( z ) | dm(z),  
$$
where $ \d_j \le 2^{-j} . $ A simple iteration proves the Lemma.
\endproof
\subsubsection*{The spaces $L^1_E $ and $H^1_E $. } 
Recall   $$L^1_E = \{ f \in L^1 ( \bT) : \hat f (k) =0 \text{ for } k
\notin E\} \quad\text{and} \quad H^1_E =L^1_E \cap H^1 . $$  
\paragraph{Embedding $L^1_E $ into $L^1( \bT^\bN ) . $}

We next define an embedding 
$$ T: L^1_E \to L^1 (\bT^\bN) $$
which maps  $H^1_E =L^1_E \cap H^1$ to the space of Hardy martingales
$H^1 (\bT^\bN).$

We define the operator by mapping the monomials $ \{z^m , m \in E \} , $ 
into $ L^1 ( \bT^\bN ) . $ Recall that for $ m \in E$ 
there exists a unique set of integers $ k_j \in ( - a_j , a_j ) $ 
so that 
$$ m = \sum k_j n _j . $$
Hence mapping the monomials  $ \{ z^m , m \in E \} , $  in $ L^1 ( \bT ) $ 
to the monomials $ \{ \prod w_j^{k_j} ,   k_j \in ( - a_j , a_j )\}  , $
gives a well    defined operator on 
$ {\rm span} \{ z^m  :  m \in E \} , $
   
$$ T :  z^m 
\to \prod w_j^{k_j} , \quad\quad m =   \sum k_j n _j . $$
Fix $ f \in L^1_E  .$  
To exhibit the martingale structure of $T (f ) $ we rewrite $T$
as follows.  
Fix
 $ w  \in \bT ^\bN $ 
and  $m \in \bN .$
Then 
write
$$
A(  k_1,  \cdots , k_{m-1},  k )  =   
\hat f ( k_1 n_1 + \cdots +k_{m-1}n_{m-1} + k n_m), $$
and 
form the Fourier coefficients

$$ a_m (k) 
 = 
\sum_{k_{m-1} = - a_{m-1}   }^{a_{m-1}}\cdots \sum_{k_1 = - a_1   }^{a_{1}} 
A(  k_1,  \cdots , k_{m-1},  k )
 w_1^{k_1}\cdots w_{m-1}^{k_{m-1}} .
$$
Thus $ a_m (k) =  a_m (k;  w_1 ,\cdots, w_{m-1}  ) .$ Then put
$$  d_m ( w)=  \sum_{k = - a_m, \, k  \ne 0    }^{a_{m}} 
  a_m (k)
w_{m}^{k} .
$$
We have  
$$T(f) (w) = \sum_{m = 1 }^\infty d_m (w) .$$
\begin{prop} \label{isoT}
For 
 $ f \in L^1_{E} $
$$   c \|f\|_{L^1( \bT ) }  \le   \| T(f) \|_{L^1  } \le  C  \|f\|_{L^1( \bT ) }  . $$
If  $ f \in H^1_E ,$ then $T(f) $  is a Hardy martingale.
\end{prop}
\proof
By inspection,  the  following properties of  $ d_m ( w) $ hold. 
First   $d_m ( w) = d_m (  w _ 1 , \dots , w_{m}),  $ second 
for fixed $ w _ 1 , \dots , w_{m-1} \in \bT $
               $$ \int_\bT d_m (  w _ 1 , \dots , w_{m-1},w_{m})d m ( w_m) = 0 , $$
and third,  if $ f \in H^1 _E ,$ then  
               $$ w_m \to  d_m (  w _ 1 , \dots , w_{m-1},w_{m}) $$
defines an analytic polynomial, hence an  element in $ H^1 _0 . $
In summary  $$ F_n
( w) = \sum_{m=0}^n d_m ( w), \quad\quad n \in \bN,  $$
gives  a Hardy martingale.
The theorem of Meyer ( see \cite{meyer} , \cite{bonami} ) asserts that for  $ f \in L^1_{E_n}( \bT )$
$$   c \|f\|_{L^1( \bT ) }   \le   \| F_n \|_{L^1 } \le  C  \|f\|_{L^1( \bT ) }  . $$
\endproof
\subsubsection*{The transfer operator  $J.$}
Define the operator  $J : L^1 ( \bT ) \to L^1 (  \bT^{\bN } ) $ by putting
$ J g= T (K*g) .$
Clearly $J$ 
is bounded 
since convolution by $K$ is a norm one operator on $ L^1 ( \bT ) $
with range on $ L^1_E ,$ and $T$ is bounded on $L^1 _E $ by
Proposition \ref{isoT}.  By Propositon \ref{small1} and Proposition \ref{isoT} 
for $ h \in L^1 ( \Sigma ) ,$
\begin{equation}\label{14sep125} \| h \|_{L^1  ( \bT ) } \le  2\| K* h
  \|_{L^1  ( \bT ) } \le 4 \| T(K*h) \|_{L^1} . \end{equation}  
hence 
$$ \frac14  \| h \|_{L^1 ( \bT ) }    \le \| J h \|_{L^1} \le 4  \| h
\|_{L^1  ( \bT )  }
, \quad\quad h \in  L^1 ( \Sigma ) . $$
Moreover by   Proposition \ref{isoT} we have the inclusion
\begin{equation}\label{14sep126}  J (  H^1 ( \bT )) \sbe  H^1 ( \bT^\bN  ). \end{equation}  
Hence by \eqref{14sep125} and  \eqref{14sep126} we proved that   $J$ satisfies \eqref{13sep121} and \eqref{13sep123}. 

\paragraph{Conditional expectation $ \bE _\cD  .$} 
Next we prove that $ \bE _\cD $ is a small perturbation of the
identity on  $J L^1 ( \Sigma ) . $
Define the non negative kernel  
$$ B (w,z) = \prod_{k = 1 }^\infty  ( 1 +  \s_k( w) \s_k ( z )) ,
\quad  w  , z  \in \bT^\bN , $$
where $ \s_k ( w)  = \s ( w_k ) . $
Conditional expectation $ \bE _\cD $ is an integral 
operator with
kernel $B ,$ 
$$   \bE _\cD  F( w ) =    \bE_z  ( B (w,z)  G(z)) , \quad \quad   F \in L^1 ( \bT^\bN ).$$ 
\paragraph{The kernels for  $J \bE_\Sigma $
 and 
$  \bE _\cD J \bE_\Sigma .$  }
Define  
the kernel 
$$A (w,\zeta) = \prod_{k = 1 }^\infty ( 1 +  s_k( w_k) \s ( \zeta^{  n_k
})), 
\quad \quad  w  \in \bT^\bN , \, \zeta \in \bT .  $$
 Then 
\begin{equation}\label{14sep127} 
J \bE_\Sigma g ( w) = \int_\bT  A (w,  \zeta)  g( \zeta  )
d m ( \zeta) , \quad  w \in \bT^\bN. \end{equation}  
Let 
$$G (w,\zeta)=
\prod_{k = 1 }^\infty \int_\bT  (1+ \s(z)\s(w_k))
( 1 +  s_k( z) \s ( \zeta ^{n_k} ))dm(z) ,\quad \quad    w \in \bT^\bN
, \zeta \in \bT . $$
Then \begin{equation}\label{14sep128} 
\bE_\cD J \bE_\Sigma g( w) = \int_\bT G (w,  \zeta)  g( \zeta )
d m ( \zeta ) 
\quad \quad  w \in \bT^\bN .  \end{equation}  
The integrals appearing in the factors of the kernel $G(w,e^{i\psi})$
may be evaluated as follows
$$
\int_\bT  (1+ \s(z)\s(w_k))
( 1 +  s_k( z) \s ( \zeta^{  n_k } ))dm(z) = 1+ \g_k\s(w_k)
\s ( \zeta ^{ n_k } ), $$
where 
$\g_k = \int_\bT \s(z)s_k( z) dm(z) . $ 
\begin{prop} Let $ g \in L^1 ( \bT ) . $ Put $ G = J \bE_\Sigma g . $
Then 
$$ \|\bE_\cD G - G \| _{L^1 } \le \e  \| g \| _{L^1 } . $$
\end{prop}\label{small2}
\proof In view of the integral representations \eqref{14sep127} and  \eqref{14sep128} 
it suffices to prove that 
$$ \sup_{\zeta \in \bT } \bE_w ( | A ( w, \zeta ) -  G( w, \zeta )|
) \le \e . $$
To this end fix $ \zeta \in \bT $ put 
$\t_k = \s ( \zeta ^{n_k } ) .$ Define
$  A_o =  G_o =1 , $ and for $j \in \bN , $

$$A_j  (w) = \prod_{k = 1 }^j ( 1 +  s_k( w_k) \t_k)),
\quad \quad G_j  (w) = \prod_{k = 1 }^j ( 1 + \g_k \s( w_k) \t_k)),  \quad\quad  w  \in \bT^\bN .$$
Rewrite the difference 
 $ A_j (w) - G_j (w) $ as 
$$
 A_{j-1} (w) \t_j (s_j( w_j) - \s ( w_j )) +     (A_{j-1} (w) - G_{j-1} (w)) \t_j ( 1 + \g_j\s(w_j)) .
$$
The second term coincides with 
$$
(A_{j-1} (w) - G_{j-1} (w)) \t_j ( 1 + \s(w_j)) -
  (A_{j-1} (w) - G_{j-1} (w)) \t_j ( 1 - \g_j)\s(w_j)).
$$
Hence 
$$ \bE | A_j  - G_j | \le ( 1 + \e _j )   \bE | A_{j-1}  - 
G_{j-1} | + 2 \e _j \bE | A_j | , $$
where we put  $\int_\bT  | s_j  - \s  | dm =\e_j . $
Since $  G_j >0 , A_j > 0 $ 
$$  \bE | A_j | = 1  , \quad \quad    \bE | G_j | = 1 , $$
Hence iterating gives 
$\bE | A_n  - G_n | \le C \sum _{j = 1}^n \e_j . $

\endproof
\paragraph{In Summary:} We proved that the linear operator
$J : L^1 ( \bT ) \to L^1 (  \bT^{\bN } ) $
defined by $ J g= T (K*g) $
 is bounded and satisfies the conditions
\eqref{13sep121}--\eqref{13sep123}. 

\bibliographystyle{abbrv}
\bibliography{hardymartingales}

\begin{thebibliography}{10}

\bibitem{bonami}
A.~Bonami.
\newblock \'{E}tude des coefficients de {F}ourier des fonctions de
  {$L^{p}(G)$}.
\newblock {\em Ann. Inst. Fourier (Grenoble)}, 20(fasc. 2):335--402 (1971),
  1970.

\bibitem{b1}
J.~Bourgain.
\newblock Embedding {$L^{1}$} in {$L^{1}/H^{1}$}.
\newblock {\em Trans. Amer. Math. Soc.}, 278(2):689--702, 1983.

\bibitem{MR827292}
J.~Bourgain.
\newblock Martingale transforms and geometry of {B}anach spaces.
\newblock In {\em Israel seminar on geometrical aspects of functional analysis
  (1983/84)}, pages XIV, 16. Tel Aviv Univ., Tel Aviv, 1984.

\bibitem{durr}
R.~Durrett.
\newblock {\em Brownian motion and martingales in analysis}.
\newblock Wadsworth Mathematics Series. Wadsworth International Group, Belmont,
  CA, 1984.

\bibitem{gar2}
D.~J.~H. Garling.
\newblock On martingales with values in a complex {B}anach space.
\newblock {\em Math. Proc. Cambridge Philos. Soc.}, 104(2):399--406, 1988.

\bibitem{g1}
D.~J.~H. Garling.
\newblock Hardy martingales and the unconditional convergence of martingales.
\newblock {\em Bull. London Math. Soc.}, 23(2):190--192, 1991.

\bibitem{jg81}
J.~B. Garnett.
\newblock {\em Bounded analytic functions}, volume~96 of {\em Pure and Applied
  Mathematics}.
\newblock Academic Press Inc. [Harcourt Brace Jovanovich Publishers], New York,
  1981.

\bibitem{sia}
A.~M. Garsia.
\newblock {\em Martingale inequalities: {S}eminar notes on recent progress}.
\newblock W. A. Benjamin, Inc., Reading, Mass.-London-Amsterdam, 1973.
\newblock Mathematics Lecture Notes Series.

\bibitem{jopfxm}
P.~W. Jones and P.~F.~X. M{\"u}ller.
\newblock Conditioned {B}rownian motion and multipliers into {${\rm
  SL}^\infty$}.
\newblock {\em Geom. Funct. Anal.}, 14(2):319--379, 2004.

\bibitem{ping}
D.~L{\'e}pingle.
\newblock Une in\'egalit\'e de martingales.
\newblock In {\em S\'eminaire de {P}robabilit\'es, {XII} ({U}niv. {S}trasbourg,
  {S}trasbourg, 1976/1977)}, volume 649 of {\em Lecture Notes in Math.}, pages
  134--137. Springer, Berlin, 1978.

\bibitem{meyer}
Y.~Meyer.
\newblock Endomorphismes des id\'eaux ferm\'es de {$L^{1}\,(G)$}, classes de
  {H}ardy et s\'eries de {F}ourier lacunaires.
\newblock {\em Ann. Sci. \'Ecole Norm. Sup. (4)}, 1:499--580, 1968.

\bibitem{m}
P.~F.~X. M{\"u}ller.
\newblock {\em Isomorphisms between {$H^ 1$} spaces}, volume~66 of {\em
  Mathematical Monographs (New Series)]}.
\newblock Birkh\"auser Verlag, Basel, 2005.

\bibitem{pfxm12}
P.~F.~X. M{\"u}ller.
\newblock A decomposition for {H}ardy martingales.
\newblock {\em Indiana. Univ. Math. J.}, pages 1--15, 2012.

\bibitem{var1}
N.~T. Varopoulos.
\newblock The {H}elson-{S}zeg{\H o} theorem and {$A_{p}$}-functions for
  {B}rownian motion and several variables.
\newblock {\em J. Funct. Anal.}, 39(1):85--121, 1980.

\end{thebibliography}
Department of Mathematics\\
J. Kepler Universit\"at Linz\\
A-4040 Linz\\
pfxm@bayou.uni-linz.ac.at

\end{document}